
\magnification=1000
\font\bfa=cmbx10 scaled 1200
\font\rmc=cmr17

\baselineskip= 16pt plus 0.6pt
\hsize=14.5truecm
\hoffset=1.5truecm
\vsize=8.6truein
\voffset=0.1truein
\def\acapo{\vskip 1pt \noindent}
\def\bib#1{\vskip 1pt\noindent\item{[#1]~\ }}

\def\Theor#1#2{\medskip \noindent{{\bf Theorem~#1.}}{#2}}
\def\Prop#1#2{\medskip \noindent{{\bf Proposition~#1.}}{#2}}
\def\Cor#1#2{\medskip \noindent{{\bf Corollary~#1.}}{#2}}
\def\Lemma#1#2{\medskip \noindent{{\bf Lemma~#1.}}{#2}}
\def\Def#1{\smallskip \noindent{{\bf Definition~#1.}}}
\def\Rem#1{\smallskip \noindent{{\it Remark.~#1.}}}
\def\Ex#1{\smallskip \noindent{{\bf Example~#1.}}}
\def\Parag#1#2{\bigskip\bigskip \noindent {\bf #1.~#2} 
\smallskip\acapo}
\def\Proof{\smallskip \noindent {{\it Proof.~}}}

\def\ninesize{scaled 1000}
\font\ninerm=cmr9 \ninesize \font\sixrm=cmr6 \ninesize 
\font\ninei=cmmi9 \ninesize \font\sixi=cmmi6 \ninesize 
\font\ninesy=cmsy9 \ninesize \font\sixsy=cmsy6 \ninesize 
\font\ninebf=cmbx9 \ninesize 
\font\nineit=cmti9 \ninesize \font\ninesl=cmsl9 \ninesize 
\skewchar\ninei='177 
\skewchar\sixi='177 \skewchar\ninesy='60 \skewchar\sixsy='60 
\def\ninepoint{\def\rm{\fam0\ninerm}
\textfont0=\ninerm \scriptfont0=\sixrm \scriptscriptfont0=\fiverm
\textfont1=\ninei \scriptfont1=\sixi \scriptscriptfont1=\fivei
\textfont2=\ninesy \scriptfont2=\sixsy \scriptscriptfont2=\fivesy
\textfont\itfam=\ninei \def\it{\fam\itfam\nineit}
\def\sl{\fam\slfam\ninesl}%
\textfont\bffam=\ninebf \def\bf{\fam\bffam\ninebf}\rm} 
\def\ds{\displaystyle} 
\def\vuoto{\emptyset}
\def\oo{\infty}

\def\pr{^{\prime }}
\def\per{\!\cdot\!}
\def\cecc#1{\mathrel{\smash{\mathop \bullet \limits _{#1}}}}

\def\Somdir#1#2{\mathrel{\smash{\mathop{\mathop \bigoplus \limits _{#1}}
\limits^{#2}}}}

\def\ldf{\leaders\hrule\hfill}
\def\riga{\kern -1.5ex\lower-.5ex\hbox to .7 truecm{\ldf}\kern -1.5ex}
\def\ldotsa{\quad\lower-.5ex\hbox {$\ldots$}\quad}
\def\mapor#1{\smash{\mathop{\longrightarrow}\limits^{#1}}}
\def\mapver#1{\Big\downarrow\rlap{$\vcenter{\hbox{$\scriptstyle#1$}}$}}

\def\solose{\Rightarrow}
\def\F{{\cal F}}\def\X{{\cal X}}\def\Y{{\cal Y}}
\def\C{{\cal C}}\def\E{{\cal E}}
\def\L{{\cal L}}
\def\O{{\cal O}}\def\M{{\cal M}}
\def\A{{\cal A}}\def\SS{{\cal S}}
\def\Hom{\hbox{\rm Hom}}\def\Ext{\hbox{\rm Ext}}\def\EXT{{\cal E}xt}
\def\HOM{{\cal H}om}

\def\dual{^{\vee}}

\def\bar#1{\,\overline{#1}\,}
\def\de{\partial} 
\def\sca{\sqcap\hskip-6.5pt\sqcup}
\def\finedim{\acapo\hfill{$\sca$}\acapo\noindent}
\mathsurround=1pt
\font\tenmsx=msbm6 scaled\magstep3 \font\sevenmsx=msbm7 \font\fivemsx=msbm5
\newfam\msxfam
  \textfont\msxfam=\tenmsx \scriptfont\msxfam=\sevenmsx
\scriptscriptfont\msxfam=\fivemsx
\font\tenmsy=msbm10 \font\sevenmsy=msbm7 \font\fivemsy=msbm5
\newfam\msyfam
  \textfont\msyfam=\tenmsy \scriptfont\msyfam=\sevenmsy
\scriptscriptfont\msyfam=\fivemsy
\font\teneuf=eufm10 \font\seveneuf=eufm7 \font\fiveeuf=eufm5
\newfam\euffam
  \textfont\euffam=\teneuf \scriptfont\euffam=\seveneuf
\scriptscriptfont\euffam=\fiveeuf
\def\Bbk{\fam\msxfam \tenmsx}

\def\com{{{\Bbk C}}}\def\raz{{{\Bbk Q}}}
\def\re{{{\Bbk R}}}\def\pro{{{\Bbk P}}}

\def\interi{{\Bbk Z}}
\def\zduer{(\interi/2)^r}\def\zdue#1{(\interi/2)^{#1}}
\overfullrule=0pt 
\bigskip\acapo
\centerline{\rmc On the moduli space of diffeomorphic algebraic surfaces}
\bigskip\bigskip\bigskip\acapo
\centerline{{\bfa Marco Manetti}}
\bigskip
\centerline{{\it Scuola Normale Superiore, Pisa, Italy}}
\bigskip
\medskip\acapo
{\ninepoint 
{\bf Abstract.} In this paper we show that the number of deformation 
types of complex structures on a fixed smooth oriented 
four-manifold can be arbitrarily large.
The examples that we consider in this paper are locally simple abelian 
covers of rational surfaces.
The proof involves the algebraic description of rational blow down, 
classical Brill-Noether theory and deformation theory of normal flat 
abelian covers.}\acapo 

\footnote{}
{\noindent
The author in member of 
GNSAGA of CNR and participates to project AGE, 
Algebraic Geometry in Europe}

\bigskip\medskip
\Parag{0}{Introduction}
One of the main problems concerning the differential topology of 
algebraic surfaces leaving unsolved by the ``Seiberg-Witten revolution''
was to determine whether the differential type of a compact complex 
surface determines the deformation type. Two compact complex manifolds 
have the same deformation type  if they are fibres of a 
proper smooth family over a connected base space (cf. [FrMo1]).
Varieties of the same deformation type are also called deformation 
equivalent.\acapo 
Here we restrict our attention to minimal surfaces of general type; in 
this case it is 
very useful to interpret the above question in terms of the moduli space 
of surfaces of general type $\M$. 
The space $\M$ is a disjoint countable union of 
quasiprojective varieties [Gi], the points of $\M$ correspond to the 
isomorphism classes of minimal surfaces of general type and two surfaces  
belong to the same connected component if and only if they are 
deformation equivalent (cf. [Gi], [Ma4], [Ca2]).\acapo
The main result of this paper is the following\smallskip
\Theor{A}{} {\it Let $X$ be a smooth oriented compact four-manifold and 
let $\M_X$ be the moduli space of minimal surfaces of general type 
orientedly diffeomorphic to $X$.\acapo
Then for every $k>0$ there exists $X$ as above with first Betti number 
$b_1(X)=0$ such that 
$\M_X$ has at least $k$ connected components.}\smallskip\acapo 
By the above remark this result gives a strong negative answer to the 
{\bf def=diff?} problem ([Ca5,\S 7], [FrMo1, Speculation of \S 3], [Ty], 
[Do]).
Note that for compact complex surfaces with $b_1=0$ and Kodaira dimension 
$<2$ the differential type determines the deformation type [Fri].
\acapo 
We recall that for an algebraic surface, $K^2$ and $\chi$ are topological 
invariants; by Gieseker theorem $\M_X$ is a quasiprojective variety and 
then it has a finite number of components.\acapo
Analogs of theorem A for the moduli space of homeomorphic complex 
surfaces already exist in [Ca4], [Ma2], [Ma3] and rely over the very 
simple algebro-geometric criterion for homeomorphism given by Freedman's 
theorem ([Fre, 1.5], [Ca1, 4.4], [Ma4, V.4]).\acapo
In view of some recent complexity results (especially [Ma3], [Ch]) and of 
the 
algebraic dependence of Donaldson and Seiberg-Witten invariants, theorem 
A is not completely surprising. 
The starting idea of this paper is the 
more subtle end effective evidence to 
theorem A given by the existence of Koll\'ar - Shepherd-Barron - 
Alexeev compactification of the moduli space of surfaces of general type
([KS], [Vie], [Ale]).
For given positive integers $a,b$ let $\M_{a,b}$ be the moduli 
space of minimal surfaces of general type $S$ 
with numerical invariants $K^2_S=a$, $\chi(\O_S)=b$; it is then possible
to embed $\M_{a,b}$ into a complete variety $\bar{\M_{a,b}^{sm}}$ 
which is a coarse 
moduli space for smoothable stable surfaces of general type [Vie, 8.39], 
[Ale], 
with numerical invariants $K^2=a$, $\chi=b$.\acapo 
As noted in [KS, 5.12] it may happen that nonhomeomorphic smooth surfaces 
belong to the same connected component of 
$\bar{\M^{sm}}=\coprod\bar{\M_{a,b}^{sm}}$; 
the substantial 
reason of this ``pathology'' is the existence of certain normal 
semilogcanonical surface singularities, each one of which  
admits at least two $\raz$-Gorenstein 
smoothings with nonhomemorphic Milnor fibres.\acapo
One can avoid this phenomenon by considering the moduli space 
$\M^T\subset\bar{\M^{sm}}$ of surfaces with at most quotient singularities.
The smoothability condition implies that only some special quotient 
singularities, called of class $T$ (Definition 1.1), 
can appear in the surfaces represented by points of $\M^T$.\acapo
By a well known result [EV] $\M^T$ is open in $\bar{\M^{sm}}$
and it is possible to prove that  
smooth minimal surfaces belonging to the same connected component of 
$\M^T$ are diffeomorphic.
Because of the extreme wilderness of 
the moduli space of surfaces of general type 
it is natural to suspect that the natural map 
$\pi_0(\M)\to\pi_0(\M^T)$ is not injective.\acapo
In practice one observe that the compactified moduli space 
$\bar{\M^{sm}}$ is a quite complicated object 
whose existence is not necessary to the proof of theorem A; 
for this reason throughout the rest of the paper we consider the simpler 
notion of deformation $T$-equivalence (Definition 1.4) instead of $\M^T$.
Roughly speaking two smooth surfaces are deformation $T$-equivalent if they 
are fibres of a proper flat map $f\colon \SS\to Y$ such that $Y$ is 
connected, 
$\SS$ is $\raz$-Gorenstein and the fibres of $f$ are normal surfaces with 
at most quotient singularities.\acapo 
The main result of section 1 (Theorem 
1.5) says that deformation T-equivalence implies smooth 
equivalence. The proof of this fact follows essentially from the 
existence of certain well-defined surgery operations on the category of 
smooth four-manifolds called ``rational blow-downs'' and already 
considered in [FS].\acapo  
The simplest singularity $X_0$ of class $T$ which is not a rational 
double point is the cone over the projectively normal rational curve of 
degree 4 in $\pro^4$; this singularity can be easily described as a 
bidouble cover $X_0\to\com^2$ with branching divisors three generic lines 
passing through $0\in\com^2$. Moreover every $\raz$-Gorenstein 
deformation of $X_0$ preserves the $\zdue{2}$-action and it is obtained 
by deforming the branching divisors (cf. 3.18).\acapo
This simple remark, together with our previous experience about abelian 
covers, motivates the use of  $\zduer$-covers in the construction of 
nontrivial examples of deformation T-equivalent surfaces.
A short but fundamental role in this paper is covered by classical 
Brill-Noether 
theory. In fact the (not yet completely developed) machinery of abelian 
covers allows to produce components of the moduli space of surfaces by 
starting from components of parameter spaces of branching divisors. The 
Brill-Noether theory gives examples of disconnected spaces $M^0_{k,n}$ 
(see \S 2) of branching divisors of the same topological type.\acapo
The detailed construction of our examples is too complicated to be 
explained in 
this introduction; however we can give here a coarse idea 
of their structure. 
Given a polarization $H$ over $Q=\pro^1\times\pro^1$,
a zerodimensional 
subscheme ( Cluster ) $\xi\subset Q$ with ideal sheaf $I_\xi$ 
is called  special if $I_\xi(H)$ is generated by global sections;
we first take the blow up $S\to Q$ along a reduced special subscheme $\xi$ 
of 
finite length $n<H^2$, then we consider  
a suitable abelian covers with group $\zdue{r}$ of $S$; in the 
construction of the cover we must take as branching divisors 
$D_\sigma\subset S$ the exceptional curves of $S\to Q$, the strict 
transform of some sections of $I_\xi(H)$ and some other very ample 
divisors. For an appropriate choice of $H,n$ and $D_\sigma$ we get 
the desired examples.\acapo 
Most part of the complexity of this construction is introduced for 
technical reasons (because the ramification of an abelian cover can 
be viewed as an obstruction to deformations and degenerations).
Together with the working construction we also give (Example 
3.20) some very simple examples of pairs of 
diffeomorphic surfaces that, at least in some case, 
we conjecture of different  deformation type.\acapo   
The paper consists of 5 sections. The first four are devoted to some 
preparatory material, most of which we consider of independent interest 
and susceptible  of future applications (cf. example 3.15). Section 5 is 
completely devoted to the construction of explicit examples of 
arbitrarily large sets of deformation T-equivalent surfaces belonging to 
different connected components of the moduli space $\M$.\acapo   
\medskip\acapo
{\it Notation.} All the varieties and schemes that we consider in this 
paper are defined over the field $\com$. For every variety $X$ we denote by 
$\Omega^1_X$, $\theta_X=\HOM(\Omega^1_X,\O_X)$ the cotangent and tangent 
sheaves respectively.\acapo
A cluster on a variety $X$ is a zero-dimensional subscheme $\xi\subset 
X$, given a cluster $\xi$ we denote by $I_\xi\subset \O_X$ its ideal 
sheaf a for a coherent sheaf $\F$ we sometimes write 
$\F-\xi=I_\xi\otimes\F$. We shall say that a cluster 
$\xi$ is $\F$-special if $\F-\xi$ is generated by global sections.\acapo
In the paper we shall use freely the following results of commutative 
algebra (cf. [Mat, \S 23]).\acapo
A) Let $G$ be a finite group acting on a local noetherian $\com$-algebra 
$A$ and let $A^G\subset A$ be the invariant subalgebra. If $A$ is normal 
(resp. Cohen-Macaulay) then $A^G$ is normal (resp. Cohen-Macaulay).\acapo
B) Let $(A,m)\to (B, n)$ be a flat morphism of local noetherian 
rings. Then $B$ is Cohen-Macaulay (resp. Gorenstein) if and only if 
$A$ and $B/mB$ are both Cohen-Macaulay (resp.: Gorenstein).
If $A$ and $B/mB$ are reduced (resp.: normal) 
then $B$ is reduced (resp.: normal).\acapo

\Parag{1}{Deformation $T$-equivalence of algebraic surfaces.}
In this section we introduce the notion of deformation $T$-equivalence of 
algebraic surfaces. It will be clear from the definition that deformation 
equivalence implies deformation $T$-equivalence, while the main result of 
this section will be the proof that deformation $T$-equivalence implies 
smooth equivalence.\acapo
\Def{1.1}([KS, 3.7]) A normal surface singularity is of {\it class $T$} 
if it is a quotient singularity and admits a one-parameter 
$\raz$-Gorenstein smoothing.\acapo
We recall that a normal  complex space is $\raz$-Gorenstein if it is 
Cohen-Macaulay and a multiple of the canonical divisor is Cartier.\acapo
A one-parameter $\raz$-Gorenstein smoothing of a normal singularity 
$(X_0,0)$ is a deformation $(X,0)\to (\com,0)$ of $(X_0,0)$ such that 
there exists a Stein representative  (as in [Lo,2.8]) 
$X\to \Delta\subset \com$ which is 
$\raz$-Gorenstein and $X_t$ is smooth for every $t\not=0$.\acapo
The singularities of class $T$ and their smoothings are well understood 
(cf. e.g. [KS], [LW], [Wa2], [Ma7]). Before recalling the classification 
we recall a standard notation concerning cyclic singularities.\acapo
If $p>0,a,b$ are integers without common irreducible factors we shall 
call {\it cyclic singularity of type $\ds{1\over p}(a,b)$} the quotient 
singularity $\com^2/\mu_p$, where the group $\mu_p=\{\xi\in\com|\, 
\xi^p=1\}$ acts on $\com^2$ by the diagonal action 
$\xi(u,v)=(\xi^au,\xi^bv)$.
It is well known (see e.g. [BPV, III.5]) that every cyclic quotient 
singularity is isomorphic to a singularity of type $\ds{1\over p}(1,q)$ 
with $p,q$ relatively prime.\acapo 
\Prop{1.2}{} {\it The singularities of class $T$ are the following:\acapo
\item{i)} Smooth points.\acapo
\item{ii)} Rational double points.\acapo
\item{iii)} Cyclic singularities of type $\ds{1\over dn^2}(1,dna-1)$ for 
$a,d,n>0$ and $a,n$ relatively prime.\acapo}\acapo
\Proof This is a well known result, for a proof see e.g. [KS, 3.11], 
[Ma7].\finedim 
\Def{1.3}{} A surface of class $T$ is a normal algebraic surfaces with at 
most singularities of class $T$.\acapo 

\Def{1.4}{} The deformation $T$-equivalence is the equivalence relation 
in the set of isomorphism classes of surfaces of class $T$ generated by 
the following relation $\sim$:\acapo
Given two surfaces of class $T$, $S_1,S_2$ we set $S_1\sim S_2$ if they 
are fibres of a proper flat analytic family $f\colon\SS\to C$ such that 
$C$ is a smooth 
irreducible curve, every fibre of $f$ is of class $T$  
and $\SS$ is $\raz$-Gorenstein.\smallskip\acapo  
It is useful to point out that 
the property of being $\raz$-Gorenstein for a one-parameter normal flat 
family 
is stable under base change and then, in the notation of 1.4, 
if $B\to C$ is a nonconstant morphism of smooth curves then 
$\SS\times_CB$ is still $\raz$-Gorenstein.\acapo   

\Theor{1.5}{} {\it If two smooth surfaces $S_1,S_2$ are deformation 
$T$-equivalent then there exists an orientation preserving diffeomorphism 
$S_1\simeq S_2$.}\acapo
Before proving 1.5 we need to recall some results about diffeomorphism of 
lens 
spaces and classification of $\raz$-Gorenstein deformations of quotient 
singularities of class $T$.\acapo 
For every oriented smooth manifold (possibly with boundary) 
$X$ we consider the following groups:\acapo  
\item{1)} $Diff(X)$ the group of diffeomorphism $X\to X$.\acapo
\item{2)} $Diff^+(X)\subset Diff(X)$ the subgroup of orientation 
preserving diffeomorphism.\acapo
\item{3)} $Diff_0^+(X)\subset Diff^+(X)$ the subgroup of orientation 
preserving diffeomorphism which are isotopic to the identity.\acapo
Let $p,q$ be relatively prime positive integers and let 
$S^3=\{(x_1,x_2)\in\com^2|\, |x_1|^2+|x_2|^2=1\}$ be the three-dimensional 
sphere.
The lens space $L(p,q)$ is by definition the quotient $S^3/\mu_p$, where 
the action is given by $\xi(x_1,x_2)=(\xi x_1, \xi^q x_2)$. The action is 
free and 
orientation preserving;  
therefore $L(p,q)$ has a natural structure of oriented 3-manifold.\acapo 
Let $\tau\colon\com^2\to\com^2$ be the complex conjugation 
$\tau(x_1,x_2)=(\bar{x_1},\bar{x_2})$, 
we have $\tau\xi=\xi^{-1}\tau$ for every $\xi\in\mu_p$ 
and therefore $\tau$ factors to an orientation preserving diffeomorphism 
$\tau\in Diff^+(L(p,q))$.\acapo
Let $\sigma\colon \com^2\to\com^2$ be the involution 
defined by $\sigma(x_1,x_2)=(x_2,x_1)$. If $q^2\equiv 1$ mod($p$) then 
$\sigma\xi=\xi^q\sigma$ and therefore $\sigma$ induces a diffeomorphism of 
$\sigma\in Diff^+(L(p,q))$. 
It is notationally convenient to define $\sigma\in Diff^+(L(p,q))$ as the 
identity if $q^2\not\equiv 1$ mod($p$).\acapo
There are other ``trivial'' diffeomorphisms of lens spaces; these are 
given by choosing a pair $(\alpha,\beta)\in S^1\times S^1$ and passing to 
the quotient the diffeomorphism 
$\varrho_{(\alpha,\beta)}(x_1,x_2)=(\alpha x_1, \beta x_2)$. 
It is immediate to observe that every diffeomorphism 
$\varrho_{(\alpha,\beta)} $ 
is isotopic to the identity.
A simple computation also shows that if $q\equiv 1$ mod($p$) (resp.: 
$q\equiv -1$ mod($p$)) then $\sigma$ (resp.: $\sigma\tau$) is conjugated 
to $\varrho_{(1,-1)}$ and then it is isotopic to the identity.\acapo
The result we need is 
\Prop{1.6}{(Bonahon)} {\it The group $\ds{Diff^+(L(p,q))\over 
Diff_0^+(L(p,q))}$ is 
generated by $\tau$ and $\sigma$.}\acapo
\Proof This is exactly Proposition 2 in [Bon]. Actually a little 
calculation, left as exercise, is needed because Bonahon consider the 
lens space as a union of two solid tori $V_1,V_2\simeq S^1\times D^2$ 
and define the 
diffeomorphisms $\sigma,\tau$ in terms  of toroidal coordinates on 
$V_1,V_2$.\finedim  
Note that if $p=dn^2$, $q=dna-1$ with $(a,n)=1$ then $q^2\equiv 1$ 
mod($p$) if and only if $a=1$, $n\le 2$.\acapo
For reader convenience we recall the notion of link of an isolated 
singularity and of Milnor fibre of a smoothing. For more details and 
proofs we refer to [Mi], [Lo] and [Wa1].
For simplicity of exposition we consider here 
only the case of normal surface singularities.\acapo 
Let $(X_0,0)$ be a normal surface singularity (hence isolated, 
irreducible and Cohen-Macaulay) and let $i_0\colon (X_0,0)\to (\com^N,0)$ 
be a closed embedding. A basic result [Mi, 2.9] asserts that there exists 
a positive real number $r$ such that for every $0<r\pr\le r$ the sphere 
$S_{r\pr}=\{z\in \com^N|\, \|z\|=r\pr\}$ intersects transversally $X_0$. 
We shall call the oriented smooth compact 3-manifold $L(X_0)=X_0\cap 
S_{r\pr}$, $0<r\pr<<1$,  the {\it link} of $(X_0,0)$. It is not difficult 
to prove  that  different choices of the embedding $i_0$ and of 
$r\pr$ give isotopic links inside $X_0-\{0\}$.
More precisely it is proved in [Lo] that for every pair $d_1,d_2\colon 
X_0\to [0,+\oo[$ of real analytic functions such that 
$d_1^{-1}(0)=d_2^{-1}(0)=\{0\}$ there exists $\epsilon>0$ such that for 
every $r_1\le \epsilon$, $r_2\le \epsilon$ the loci 
$d_1^{-1}(r_1)$ and $d_2^{-1}(r_2)$ are isotopic smooth subvarieties of 
$X_0-\{0\}$.\acapo
Note that $L(p,q)$ is the link of the cyclic singularity $X_0$ of type 
$\ds{1\over p}(1,q)$; in fact it is sufficient to take as real analytic 
function $d\colon X_0\to[0,+\oo[$ the usual norm of $\com^2$ which is 
clearly $\mu_p$ invariant.\acapo
Assume now that $f\colon (X,0)\to (\com,0)$ is a one-parameter smoothing 
of 
$(X_0,0)$ and let $i\colon (X,0)\to (\com^N\times\com,0)$ be a closed 
embeddings extending $i_0$ such that $f$ is the composition of $i$ with 
the projection on the second factor. Let $r>0$ as above, for every 
$0<r\pr\le r$ there exists a $\delta>0$ such that $X_t=f^{-1}(t)$ 
intersects transversally $S_{r\pr}$ for every $t\in\com$, $|t|\le 
\delta$.\acapo 
If $0<|t|<<1$ the smooth oriented  manifold with boundary 
$F=X_t\cap \{z\in\com^N| \|z\|\le r\pr\}$ is called the {\it Milnor 
fibre} of the smoothing. Again its diffeomorphism class is independent 
from the choice of $i,r\pr$ an $t$ and by Ehresmann fibration theorem 
$\de F=L(X_0)$. If $(X_0,0)$ is a complete intersection then the 
diffeomorphism type of $F$ is also 
independent from the smoothing [Lo], in general $F$ is not uniquely determined 
by $(X_0,0)$.\acapo
By Lefschetz duality we have $H^2(F,\interi)=H_2(F,\de F,\interi)$ and 
the natural map $H_2(F,\interi)\to H_2(F,\de F,\interi)$ gives the 
intersection product $q\colon H_2(F,\interi)\times 
H_2(F,\interi)\to\interi$.\acapo  
It is also trivial to see that the Milnor fibre  
is invariant under base change and if the smoothing admits simultaneous 
resolution [Tya] then $F$ is isomorphic to a neighbourhood of the 
exceptional 
curve of the minimal resolution of $X_0$.\acapo
\Prop{1.7}{} {\it The Milnor fibre $F$ of a $\raz$-Gorenstein 
smoothing of a cyclic singularity of class $T$ is unique up to 
orientation preserving diffeomorphism. Moreover the restriction 
homomorphism $Diff^+(F)\to Diff^+(\de F)$ is surjective.}\acapo
\Proof  As $\de F$ is collared in $F$ the restriction map 
$Diff_0^+(F)\to Diff_0^+(\de F)$ is surjective. As $\de F$ is a lens space, 
according to 1.6 it is sufficient to prove that $\sigma,\tau$ lift to 
$Diff^+(F)$.\acapo
Let $d,n,a$ be fixed positive integers with $(a,n)=1$ and let $X_0$ be 
the cyclic quotient singularity of type $\ds{1\over dn^2}(1,dna-1)$. 
If $\xi\in\mu_{dn^2}$ then $\xi$ acts on $\com^2$ by a 
linear transformation of determinant $\xi^{dna}$ 
and then the subgroup $\mu_{dn}\subset\mu_{dn^2}$ acts on $\com^2$ by 
linear transformations preserving the holomorphic form $dx_1\wedge dx_2$. 
In particular 
$X_0$ can be 
described as $Y_0/\mu_n$, where $Y_0=\com^2/\mu_{dn}$ 
is the rational double point of type $A_{dn-1}$ (i.e. the cyclic 
singularity of type $\ds{1\over dn}(1,-1)$).\acapo
Setting $u=x_1^{dn}$, $v=x_2^{dn}$ and $y=x_1x_2$
we can describe $Y_0$ as the hypersurface singularity of $\com^3$ defined 
by the equation $uv-y^{dn}=0$; in this coordinates $\mu_n$ acts by 
$$\mu_n\ni\xi\colon (u,v,y)\to (\xi u,\xi^{-1}v,\xi^a y)$$
Moreover  the diffeomorphisms 
$\tau,\sigma$ of $L(Y_0)$ are induced by $\tau,\sigma\in Diff(\com^3)$
$$\tau(u,v,y)=(\bar{u},\bar{v},\bar{y}),\qquad 
\sigma(u,v,y)=(v,u,y).$$
It is known [KS, 3.17],  that every $\raz$-Gorenstein one-parameter 
deformation of $(X_0,0)$ is isomorphic to the pull-back, via a suitable  
holomorphic 
germ map $(\com,0)\to (\com^d,0)$, of the deformation 
$\pi\colon \Y/\mu_n\to\com^d$, where $t_0,...,t_{d-1}$ are coordinates 
over $\com^d$, $\Y\subset \com^3\times\com^d$ is the hypersurface of 
equation $uv-y^{dn}=\sum_{k=0}^{d-1}t_ky^{kn}$ and $\mu_n$ acts as 
$$\mu_n\ni\xi\colon (u,v,y,t_0,...,t_{d-1})\to 
(\xi u,\xi^{-1}v,\xi^a y, t_0,...,t_{d-1}).$$
Since $\com^d$ is locally irreducible, the basic theory of the 
discriminant 
tell us that all the Milnor fibres of $\raz$-Gorenstein smoothings of 
$X_0$ 
are orientedly diffeomorphic. It is therefore sufficient to prove that 
$\sigma,\tau$ extends to diffeomorphism of $Y_t/\mu_n$ for some $t\in 
\com^d$; this is immediate to check if $t=(t_0,0,..,0)$ for $t_0$ a real 
number $0<t_0<<1$.\finedim
\Ex{1.8} A simple calculation which we omit shows that the Milnor fibre 
of a $\raz$-Gorenstein smoothing of    
the cyclic singularity of type $\ds{1\over 4}(1,1)$ is the complement in 
$\pro^2$ of a tubular neighbourhood of the conic of equation $uv-y^2=0$. 
The diffeomorphism $\tau,\sigma$ are induced by the diffeomorphism of 
$\pro^2$
$$\tau(u,v,y)=(\bar{u},\bar{v},\bar{y}),\qquad 
\sigma(u,v,y)=(v,u,y).$$\smallskip\acapo
{\it Proof of 1.5.~~} We prove 1.5 by associating to every surface $S$ of 
class 
$T$ an oriented smooth 4-manifold $\hat{S}$ which is equal to $S$ when 
$S$ has no singular points; next we will prove that if $S_1\sim S_2$, 
where $\sim$ is the relation introduced in 1.4, then $\hat{S_1}$ is 
diffeomorphic to $\hat{S_2}$.\acapo
As a first step we define for every surface of class $T$, 
$\hat{S}=\hat{S\pr}$ where $S\pr\to S$ is the minimal resolution of all 
rational double points of $S$.\acapo
Let $S$ be a surface with at most cyclic singularities of class $T$ at 
points $p_1,...,p_n$. Choose closed embeddings 
$\phi_i\colon(S,p_i)\to (\com^{N_i},0)$, $r>0$ a sufficiently small real 
number and 
denote $V_i=\phi_i^{-1}(\{z|\, \|z\|< r\})$. If $F_i$ is the Milnor fibre 
of a $\raz$-Gorenstein smoothing of $(S,p_i)$ we define
$$\hat{S}=(S-\cup_{i}V_i)\cup(\cup_{i}F_i)$$
where the pasting is made by choosing for every $i$ an orientation 
preserving 
diffeomorphism $\de V_i\to\de F_i$. By proposition 1.7 $\hat{S}$ is a 
well defined smooth oriented four manifold.\acapo
By a simple topological argument it is sufficient to 
prove that if $f\colon \SS\to\Delta=\{t\in\com|\, |t|<1\}$ is a 
$\raz$-Gorenstein deformation of a surface $S_0$ of class $T$ then 
$\hat{S_0}=\hat{S_t}$ for $|t|<<1$.\acapo
By using the simultaneous resolution of rational double points 
we can assume without loss of generality that $S_0$ has at most cyclic 
singularities at points $p_1,...,p_n$.\acapo
Step 1). Assume first that $S_t$ is smooth for $t\not=0$,  
$\phi_i\colon(\SS,p_i)\to (\com^{N_i}\times \Delta,0)$  closed embeddings; 
then for 
$r>0$ sufficiently small, $V_i=S_0\cap \phi_i^{-1}(B_r)$, 
$F_i=S_t\cap \phi_i^{-1}(B_r)$, $t\not=0$. By integration of vector 
fields it is easy to construct, for $0<|t|<<1$, a diffeomorphism 
$(S_0-\cup_{i}V_i)\to (S_t-\cup_{i}F_i)$ sending $\de V_i$ into $\de F_i$
and then by proposition 1.7 
$$\hat{S_0}=(S_0-\cup_{i}V_i)\cup(\cup_{i}F_i)
=(S_t-\cup_{i}F_i)\cup(\cup_{i}F_i)=S_t.$$ 
Step 2). Assume now $S_t$ possibly singular and let 
$p\in S_0$ be a singular point; 
by using the classification theorem of $\raz$-Gorenstein deformations of 
singularities of class $T$ (see the proof of 1.7) it is immediate to find 
a $\raz$-Gorenstein deformation $\tilde{f}\colon (\tilde{\SS},p)\to 
(\Delta_t\times \Delta_s,0)$ such that $f$ is the pull-back of $\tilde{f}$ 
via the inclusion $\Delta_t\times\{0\}\subset \Delta_t\times\Delta_s$ and  
$S_{t,s}$ is smooth for $s\not=0$. Therefore $\tilde{\SS}$ gives a local 
simultaneous $\raz$-Gorenstein smoothing
of the fibres of $f$. Step 1) and a simple pasting 
argument concludes the proof.
\finedim
{\it Remark.} Our first proof of 1.5 didn't make use of proposition 1.7 
but used a more careful study of the Kuranishi families of singularities 
of class $T$; however this approach 
(i.e. assuming Bonahon's theorem) is considerably 
simpler and more elegant.\acapo  
If $V\to S$ is the minimal resolution of a surface of class $T$ then, in 
the terminology of [FS],  
$\hat{S}$ is obtained by rationally blowing down $V$.\acapo

\Parag{2}{Special clusters  on $\pro^1\times\pro^1$ and related spaces.}
Let $Q=\pro^1\times\pro^1$ be the quadric  and  
let $H=\O_Q(a,b)$ be a fixed polarization over $Q$ with $a,b\ge 3$; we 
have $H^2=2ab$.\acapo
For every pair of integers $k\ge 0$, $0\le n\le 2ab$, let 
$M_{k,n}^{a,b}$ be the  subset of $|H|^k\times Q^n$ consisting of the elements 
$(C_1,...,C_k,p_1,...,p_n)$ such that:\acapo
\item{i)} the curves $C_i$'s are smooth.\acapo
\item{ii)} For every $i\not=j$ $C_i$ intersect transversally $C_j$.\acapo
\item{iii)} $p_i\in C_j$ for every $i,j$, i.e. the points $p_i$ are 
contained in the base locus of the linear system 
generated by $C_1,...,C_k$.\acapo
\item{iv)} $p_i\not=p_j$ if $i\not=j$.\acapo
For simplicity of notation we shall write 
$M_{k,n}$ instead of $M_{k,n}^{a,b}$
whenever there is non ambiguity about the polarization.\acapo
The set $M_{k,n}$ carries a natural structure of locally closed subscheme 
of $|H|^k\times Q^n$.\acapo 
For $l\le k$, $m\le k$ we shall call natural projection $M_{k,n}\to 
M_{l,m}$ the map 
$$(C_1,...,C_k,p_1,...,p_n)\to 
(C_1,...,C_l,p_1,...,p_m);$$ 
it is clearly a regular morphism of schemes.\acapo  
If $Aut_0(Q)=PGL(2)\times PGL(2)$ denotes the group of biregular 
automorphisms of $Q$ acting trivially over $Pic(Q)$ then there exists a 
natural regular action of $ Aut_0(Q)\times \Sigma_k\times \Sigma_n$ 
over $M_{k,n}$, where $\Sigma_i$ is the symmetric group of permutations 
of $i$ elements.\acapo

\Lemma{2.1}{} {\it The scheme $M_{k,n}$ is connected, if $k\le 2$ it is 
also  smooth irreducible.}\acapo
\Proof If $k\le 1$ or $n=0$ the lemma is trivially true. By the implicit 
function theorem the natural projection $M_{2,n}\to M_{2,0}$ is an 
unramified covering, it is therefore sufficient to prove 
that there exists a fibre 
contained in a path-connected subset of $M_{2,n}$.\acapo
Let $C\in |H|$ be a fixed smooth curve, as $|H|$ cuts over $C$ a very 
ample linear system, by the general position theorem (see below) the 
fibre over $C$ of the projection $M_{2,n}\to M_{1,0}$ is connected.\acapo
If $k\ge 3$ let $V\subset M_{k,n}$ be the subset of elements  
$(C_1,...,C_k,p_1,...,p_n)$ such that every curve $C_i$ belong 
to the pencil generated by $C_1,C_2$. The natural projection $V\to 
M_{2,n}$ is a smooth map with irreducible fibres, in particular $V$ is 
also smooth irreducible. The conclusion follows from the fact that it is 
always possible to join by a path every point of $M_{k,n}$ with a point 
of $V$.\finedim 
In the above proof we have used the following version of the general 
position theorem; for a proof we refer to [ACGH] p. 112.\acapo
\medskip\acapo
{\bf General position theorem.}{} {\it Let $C\subset \pro^r$ be a smooth 
curve of degree $d$. Then 
$$I=\{(p_1,...,p_d,H)\in C^d\times(\pro^r)\dual|
\, p_i\not=p_j,\, \{p_1,...,p_d\}=H\cap C\}$$
is smooth irreducible of dimension $r$; hence
for every $s\le d$ the image of the projection of $I\to C^s$, 
$(p_1,...,p_d,H)\to (p_1,...,p_s)$ is 
irreducible too.}\medskip\acapo
Let $M^0_{k,n}\subset M_{k,n}$ be the (possibly empty) open subset of 
elements $(C_1,...,C_k,p_1,...,p_n)$ such that the base locus of the 
linear system generated by $C_1,..,C_k$ is {\bf exactly} 
$\{p_1,...,p_n\}$.\acapo 
Consider now a positive integer $0<c<\ds{1\over 2}a$ and let 
$L=\O_Q(a-c,b)$, $F=\O_Q(c,0)=H-L$.\acapo
Let $n=H\per L=b(2a-c)$ and consider for every $k\ge 1$ 
$$M_{k,L}=\{(C_1,...,C_k,p_1,...,p_n)\in M_{k,n}
|\,\O_{C_1}(\sum p_i)=\O_{C_1}(L)\}$$
Note that if $c\ge 2$ then the equality $\O_{C_1}(\sum p_i)=\O_{C_1}(L)$ 
does not imply that the divisor $p_1+...+p_n\subset C_1$ is cutted by a 
curve in the linear system $H^0(Q,L)$.\acapo
$M_{k,L}$ has a natural structure of closed subscheme of $M_{k,n}$, in 
order to show this, using the fact that $M_{k,L}$ is the fibred product of 
$M_{1,L}\to M_{1,n}$ and the natural projection $M_{k,n}\to M_{1,n}$, it 
is not restrictive to assume $k=1$.\acapo
Let $\C\mapor{\pi}M_{1,n}$ be the universal curve, $\C\subset Q\times 
M_{1,n}$, $q\colon \C\to Q$ the projection
and let $p_i\colon M_{1,n}\to \C$ be the sections of marked 
points, then, as the fibres of $\pi$ are smooth irreducible curves,  
$M_{1,L}$ is naturally isomorphic to the first Fitting subscheme of 
the coherent sheaf $\pi_*\O_{\C}(q^*L-\sum p_i)$. We consider the Fitting 
subscheme instead of the support because we want that $M_{k,L}$ represent 
the corresponding functor from the opposite category of schemes over 
$\com$ to the category of sets.

\Lemma{2.2}{} {\it Let $C_1\in |H|$ be a smooth curve, then the restriction 
map $H^0(Q,F)\to H^0(C_1,F)$ is an isomorphism. In particular if 
$(C_1,...,C_k,p_1,...,p_n)\in M_{k,L}$ then $\O_{C_i}(\sum 
p_i)=\O_{C_i}(L)$ for every $i=1,...,k$
and $M_{k,L}$ is stable under the action of 
$Aut_0(Q)\times \Sigma_k\times \Sigma_n$.}\acapo
\Proof The first part follows immediately from the vanishing of 
$H^0(Q,-L)$ and $H^1(Q,-L)$. For the second part it is not restrictive to 
assume $k=2$. Write $C_1\cap C_2=\sum p_i+D$ with $D\in|\O_{C_1}(F)|$ a 
reduced divisor, by the first part of the lemma there exists a curve 
$F_1\in|F|$ such that $D=C_1\cap F_1$. 
Since $D$ is reduced over $C_1$ 
of degree $cb$ we also have $D=C_2\cap F_1$ and 
therefore $\O_{C_2}(\sum p_i)=\O_{C_2}(H-F)=\O_{C_2}(L)$.\finedim
\Prop{2.3}{} {\it In the above notation the scheme 
$M_{k,L}$ is smooth irreducible of dimension 
$n+2(a+b)-1+(k-1)(c+1)$. Moreover the natural projections 
$M_{k,L}\to M_{1,L}\to M_{1,0}$ are smooth morphisms.}\acapo
\Proof The space $M_{1,0}$ is a Zariski open subset of $|H|$ and 
therefore it is smooth irreducible of dimension $ab+a+b$.\acapo
We have $H^1(Q,L)=H^2(Q,L-H)=0$ 
and therefore for every smooth curve $C\in |H|$ we have $g(C)=(a-1)(b-1)$ 
and 
$H^1(C,L)=0$.\acapo 
By Riemann-Roch formula
$h^0(C,L)=L\per C+1-g(C)=b(a-c)+a+b$.
Let $C\mapor{p} M_{1,0}$ be the universal curve, then $p_*q^*L$ is a 
locally free sheaf of rank $b(a-c)+a+b$. 
Now $M_{1,L}$ is an unramified covering of an 
open subset of the projectivized of $p_*q^*L$ and 
therefore it is smooth of dimension $n+2(a+b)-1$. The same proof shows 
also that the natural projection $M_{1,L}\to M_{1,0}$ is smooth and by 
the general position theorem its fibres are irreducible.\acapo
In order to conclude the proof we show that the fibres of 
the projection $M_{k,L}\mapor{\pi} M_{1,L}$ are smooth irreducible of 
dimension 
$(k-1)(c+1)$.\acapo 
First note that for every $m=(C,p_1,...,p_n)\in M_{1,L}$ 
the complete linear system $|\O_C(H-\sum p_i)|=|\O_C(F)|$ is base point 
free of dimension $c$. From this it follows that  
the fibre of $\pi$ over 
$m$ is a nonempty open subset of $\pro(H^0(H-p_1-...-p_n))^{k-1}$.\finedim
We are now able to prove the main result of this section 
\Theor{2.4}{} {\it In the notation above, for every $k\ge 3$, 
$M_{k,L}^0:=M_{k,L}\cap M_{k,n}^0$ 
is a nonempty connected component of $M_{k,n}^0$.}\acapo
\Proof 
By proposition 2.3 the space $M_{k,L}$ is smooth irreducible. 
Given generic curves $C\in |H|$, $D\in |L|$, $F_1,F_2\in |F|$ and 
setting $\{p_1,...,p_n\}=C\cap D$ we have  
by Bertini's theorem that 
$(C_1,...,C_k,p_1,...,p_n)\in M_{k,L}^0$, where $C_1,...,C_k$ are generic 
curves in the twodimensional linear system generated by 
$C, D+F_1$ and $D+F_2$.\acapo
The scheme  $M_{k,L}^0$ is the restriction of $M_{k,L}$ to
$M^0_{k,n}$ and it is therefore closed is the latter; by implicit 
function theorem it is sufficient to prove that for every 
$m=(C_1,...,C_k,p_1,...,p_n)\in M^0_{k,L}$ 
the natural map of Zariski tangent spaces $T_{m,M_{k,L}}\to 
T_{m,M_{k,n}}$ is surjective.\acapo
We have already proved that the composition 
$M_{k,L}\mapor{i} M_{k,n}\mapor{\alpha} M_{1,0}$ is smooth; we now prove 
that every morphism $\phi\colon T=Spec(\com[\epsilon]/(\epsilon^2))\to 
(M_{k,n},m)$ such that $\phi\alpha=0$ can be lifted to   
$M_{k,L}$.\acapo 
Note that for every line bundle $\L$ over $Q$, 
the line bundle $q^*\L$ is the unique extension of $\L$ to $Q\times 
T$, where $q\colon Q\times 
T\to Q$ is the projection.\acapo
The morphism $\phi$ represents $k$ smooth curves over $T$,
$C_{T,1},...,C_{T,k}$ which are divisors of 
$q^*H$ together with $n$ sections 
$p_{T,i}\colon T\to \cap C_{T,j}\subset Q\times T$ 
such that their restriction to the closed 
point of $T$ give $m$; the condition $\alpha\phi=0$ is equivalent to 
$C_{T,1}=C_1\times T$.\acapo 
For every $i=2,...,k$ let $D_i=C_1\cap C_i-\sum p_j\in |\O_{C_1}(F)|$ 
and let $s_i\in H^0(Q,F)$ be a section whose restriction to $C_1$ gives 
the divisor $D_i$. The assumption $m\in M_{k,L}^0$ implies that the 
linear system $V\subset H^0(Q,F)$ generated by $s_2,...,s_k$ has no fixed 
components and since $F^2=0$ it has no base points.\acapo
Denoting by $W\subset H^0(C_1,F)$ the isomorphic image of $V$ under the 
restriction map, we have a commutative diagram with vertical epimorphisms and 
horizontal cup products
$$\matrix{
V\otimes H^0(Q,K_Q+L)&\mapor{\beta}&H^0(Q,K_Q+C_1)\cr
\mapver{}&&\mapver{}\cr
W\otimes H^0(C_1,K_{C_1}-F)&\mapor{\mu_0}&H^0(C_1,K_{C_1})\cr}$$
The map $\beta$ is surjective; in fact it is well known (and easy to 
prove) that, denoting  $S_d\subset\com[x,y]$ the subspace of homogeneous 
polynomials of degree $d$,  if $V\subset S_c$ is a base point free (over 
$\pro^1$) vector subspace
then  $V\otimes S_{a-c}\to S_a$ is surjective for every $a\ge 2c-1$. 
Moreover it is also well known that if $A,B$ are effective divisor over 
$Q$ then the cup product $H^0(Q,A)\otimes H^0(Q,B)\to H^0(Q,A+B)$ is 
surjective; in conclusion $\beta$ and then $\mu_0$ are surjective.\acapo
The map $\mu_0$ plays an important role in the description of the 
Zariski tangent space of the Brill-Noether varieties $W^r_d(C_1)$ and 
$G^r_d(C_1)$ (see [ACGH] for the definitions) at the points $\O_{C_1}(F)$ 
and $W$ respectively.\acapo 
In particular by [ACGH, IV.4.1], the 
annihilator $Im(\mu_0)^{\perp}\subset H^1(\O_{C_1})=H^0(K_{C_1})\dual$ 
classifies the first order extension of $\O_{C_1}(F)$ where the divisors 
$D_i$ extend. As $\mu_0$ is surjective
it follows 
immediately that $\O_{C_1\times T}(q^*H-\sum p_{T,i})=\O_{C_1\times 
T}(q^* F)$ and then $\phi$ can be lifted to $M_{k,L}$.\finedim   
By simmetry the above results are still true if 
$L=\O_Q(a,b-c\pr)$ with $2c\pr<b$.
The proof of 2.4 does not work if $L=\O_Q(a-c,b-c\pr)$, $c,c\pr>0$, even 
in the case $a>>c, b>>c\pr$, because in general $V$ is not base point 
free and the map $\mu_0$ cannot be surjective.\acapo 

\Cor{2.5}{} {\it Assume $b=la$ for some integer $l\ge 2$ and $n=la(2a-c)$, 
$0<2c<a$, $k\ge 3$.\acapo 
Then $M^0_{k,n}$ contains at least two connected components which are 
stable under the action of $Aut_0(Q)\times \Sigma_k\times 
\Sigma_n$.}\acapo 
\Proof Setting $L_1=\O_Q(a-c,b)$, $L_2=\O_Q(a,b-lc)$ 
we have by 2.3 $M_{k,L_1}\cap M_{k,L_2}=\vuoto$ and then  
we can take 
$M^0_{k,L_1}$ and $M^0_{k,L_2}$
as connected components.\acapo
A simpler proof of the equality $M_{k,L_1}\cap M_{k,L_2}=\vuoto$ 
which works also for 
$k=1$ follows from 2.2. In fact if 
$(C_1,...,C_k,p_1,...,p_n)\in M_{k,L_1}\cap M_{k,L_2}$
then 
$$h^0(\O_{C_1}(H-\sum p_i))=\cases{h^0(\O_{C_1}(H-L_1))=c+1&\cr 
h^0(\O_{C_1}(H-L_2))=lc+1&\cr}.$$
~\finedim 
By a careful reading of the proofs we see that all the above results 
concerning the local structure of the spaces $M_{k,L}$ and $M_{k,n}$ at 
a point $(C_1,...,C_k,p_1,...,p_n)$ depends only on the properties of the 
linear system generated by the curves $C_1,..,C_k$. 
We can generalize these results by considering the ``closure''
$\bar{M}_{k,n}\subset |H|^k\times Q^n$ whose elements 
$(C_1,...,C_k,p_1,...,p_n)$ satisfy the conditions:\acapo
\item{i)} the curves $C_i$'s are reduced without common components.\acapo
\item{ii)} $p_i\in C_j$ for every $i,j$, i.e. the points $p_i$ are 
contained in the base locus of the linear system $V$ generated by 
$C_1,...,C_k$.\acapo
\item{iii)} $p_i\not=p_j$ if $i\not=j$.\acapo
\item{iv)} The subscheme of base points of $V$ is reduced, 
i.e. if $q$ belongs to the intersection of $C_1,..,C_k$ then 
the local equations of the curves $C_i$'s at the 
point $q$ generate the maximal ideal.\acapo
We define $\bar{M}_{k,n}^0\subset \bar{M}_{k,n}$ as the open subscheme of 
points $(C_1,...,C_k,p_1,...,p_n)$ such that 
$\{p_1,...,p_n\}=C_1\cap...\cap C_k$.\acapo
It is immediate to see, by using Bertini's theorem, that if 
$(C_1,...,C_k,p_1,..,p_n)\in \bar{M}_{k,n}$  
(resp.: $\bar{M}_{k,n}^0$) then for generic curves 
$D_1,...,D_k$ in the linear system generated by the divisors $C_i$ we have 
$(D_1,...,D_k,p_1,..,p_n)\in M_{k,n}$ (resp.: $M_{k,n}^0$).\acapo
We then define 
$\bar{M}_{k,L}$ as the set of $(C_1,...,C_k,p_1,...,p_n)\in 
\bar{M}_{k,n}$ such that $\O_{D}(\sum p_i)=\O_{D}(L)$ for some smooth 
curve $D$ in the linear system generated by $C_1,...,C_k$. By  
Lemma 2.2 if $\O_{D}(\sum p_i)=\O_{D}(L)$ for one smooth curve then the 
same is true for every smooth curve in the linear system 
$<C_1,...,C_k>$.\acapo 
As above $\bar{M}_{k,L}$ and 
$\bar{M}_{k,L}^0=\bar{M}_{k,L}\cap \bar{M}_{k,n}^0$ 
are closed subschemes  of 
$\bar{M}_{k,n}$, $\bar{M}_{k,n}^0$ respectively.\acapo
\Theor{2.6}{} {\it $\bar{M}_{k,L}$ is a smooth irreducible variety and if 
$m=(C_1,...,C_k,p_1,...,p_n)\in \bar{M}_{k,L}^0$  then 
$\bar{M}_{k,L}\to \bar{M}_{k,n}$ is a local isomorphism at $m$.}\acapo
\Proof 
It is convenient to work with the $(\com^*)^k$-principal bundles 
$F_{k,n}\to M_{k,n}$, $F_{k,L}\to M_{k,L}$, 
$\bar{F}_{k,n}\to \bar{M}_{k,n}$, $\bar{F}_{k,L}\to \bar{M}_{k,L}$, 
obtained by considering the sections of the line bundle $H$ instead of 
the divisors.\acapo
Let $m=(f_1,...,f_k,p_1,...,p_n)\in\bar{F}_{k,n}$ be a fixed point and let 
$A=(a_{ij})\in GL(k,\com)$ be a matrix, we then define 
$Am=(g_1,...,g_k,p_1,...,p_n)$ where $g_i=\sum a_{ij}f_j$. If $A$ is 
generic then $Am\in F_{k,n}$ and $A$ extends to an isomorphism of germs 
of complex spaces $A\colon (\bar{F}_{k,n},m)\to (F_{k,n},Am)$; this germ 
isomorphism clearly preserves the subscheme $\bar{F}_{k,L}$ and the 
theorem follows by the above results.\finedim
We also point out that $M_{k,n}$ is Zariski open and dense in 
$\bar{M}_{k,n}$. The motivation of our definition of $\bar{M}_{k,n}$ and 
$\bar{M}_{k,n}^0$ will be clear in the next sections.
Moreover the inclusion $M_{k,n}^0\subset\bar{M}_{k,n}^0$ induces a 
bijection between the sets $\pi_0$'s of connected components.
\acapo 
Consider now a fixed $m=(C_1,...,C_k,p_1,...,p_n)\in\bar{M}_{k,n}$ with 
$k\ge 2$ and let 
$S=S_m\mapor{\nu}Q$ be the blow up at the points $p_1,...,p_n$.
If $E_i\subset S_m$ is the $(-1)$-curve over $p_i$ and  
$C_i\pr=\nu^*C_i-\sum_j E_j$ then 
$m\in \bar{M}_{k,n}^0$ if and only if $\cap C_i\pr=\vuoto$.\acapo

\Lemma{2.7}{} {\it In the above notation 
$h^0(\theta_S)=h^2(\theta_S)=0$, $h^1(\theta_S)=2n-6$.}\acapo
\Proof The vanishing of $H^2(\theta_S)$ is a general fact which holds for 
every smooth rational surface.\acapo
Let $E_i\subset S$ be the exceptional curve over the point $p_i$, 
$1=1,..,n$. It is well known that there exists two exact sequences of 
sheaves
$$0\mapor{} 
\theta_S\mapor{}\nu^*\theta_Q\mapor{}\oplus\O_{E_i}(1)\mapor{}0$$
$$0\mapor{} \nu_*\theta_S\mapor{}\theta_Q\mapor{}\oplus 
T_{p_i,Q}\mapor{}0$$
and we also have 
the vanishing of the higher direct image sheaf $R^1\nu_*\theta_S=0$.\acapo
It is therefore sufficient to prove that the natural restriction map 
$H^0(Q,\theta_Q)\to \oplus T_{p_i,Q}$ is injective. Let $\pi_i\colon Q\to 
\pro^1$, $i=1,2$ be the natural projections; 
as $n>\max(2a,2b)$ and the curves $C_i$ have no common components, 
the points $p_i$ lie in at least three different fibres of $\pi_1$, 
$\pi_2$. 
Since $\theta_Q=\O_Q(2,0)\oplus\O_Q(0,2)$ the injectivity of the above 
restriction map is immediate.\finedim
By a well known formula we have an isomorphism 
$Pic(S)=\nu^*Pic(Q)\oplus\interi E_1\oplus...\oplus\interi E_n$; from now 
on we shall write $\O_S(r,s)=\nu^*\O_Q(r,s)$. By the Leray spectral 
sequence we get $h^i(\O_S(r,s))=h^i(\O_Q(r,s))$.\acapo
\Def{2.8}{} We shall say that a line bundle $\L=\O_S(r,s)-\sum a_iE_i$ is 
combinatorially ample if 
$a_i\ge 2$ for every $i=1,...,n$ and $r,s> \sum(a_i+1)$.\acapo
Note that if $\L$ is combinatorially ample then for every  
$D\in Pic(S)$ the line 
bundle $D+\alpha\L$ is combinatorially ample for 
$\alpha$ sufficiently large.\acapo
\Prop{2.9}{} {\it Let $\L\in Pic(S)$ be a combinatorially ample line 
bundle, then:\acapo
\item{i)} $\L$ is ample and base point free.\acapo
\item{ii)} $H^1(\L)=H^1(-\L)=0$.\acapo
\item{iii)} $H^0(\theta_S(-\L))=H^1(\theta_S(-\L))=0$.\acapo}
\Proof Consider a line bundle $\L=\O_S(r,s)-\sum a_iE_i$ 
with $a_i>0$ for every $i$.\acapo
a) If $r,s\ge \sum a_i$ then $\L$ is effective and nef. In fact 
$\L^2=2rs-\sum a_i^2> 0$, $\L\per E_i=a_i>0$; moreover if $V\subset 
H^0(\L)$ is the subspace generated by sections whose zero locus is 
contained in a finite union of fibres of the two projections 
$\pi_1,\pi_2\colon S\to \pro^1$, then it is immediate to observe that the 
possible fixed part of $V$ is contained in the exceptional divisor 
$E=\cup E_i$. In particular $\L\per C\ge 0$ for every irreducible curve on 
$S$.\acapo 
b) If $r,s>\sum a_i$ then $\L$ is ample. This follows by point a) and 
Nakai criterion.\acapo
c) If $\delta\colon X\to S$ is the blow up at a point with exceptional 
curve $E\subset X$ the same arguments of points a) and b) show that if 
$r,s\ge 3+\sum a_i$ then $\delta^*\L-E$ is ample.\acapo
d) By adjunction formula $K_S=\O_S(-2,-2)+\sum E_i$ and by point b), if 
$\L$ is combinatorially ample then $\L$, $\L-K_S$ are ample. By Kodaira 
vanishing we have $H^1(-\L)=0$, $H^1(\L)=H^1(K_S-\L)\dual=0$.\acapo
The same arguments shows that $\delta^*(\L-K_S)-2E$ is ample on $X$, in 
particular $H^1(X,K_X-\delta^*\L+E)=0$ and $\L$ is base point free on 
$S$.\acapo 
e) We have an exact sequence
$$0\mapor{} \theta_S(-\L)\mapor{}(\O_S(2,0)-\L)\oplus(\O_S(0,2)-\L)
\mapor{}\oplus\O_{E_i}(1-a_i)\mapor{}0.$$
By taking the associate cohomology exact sequence, if $\L$ is 
combinatorially ample then by the above steps we get 
$H^0(\theta_S(-\L))=H^1(\theta(-\L))=0$.\finedim
Let $\SS\mapor{\nu}\bar{M}_{k,n}\times Q$ be the blow up of the sections 
$p_i\colon \bar{M}_{k,n}\to Q$. It is easy to see, using the fact that 
the union of the images of $p_1,...,p_n$ is a locally complete 
intersection subscheme of $\bar{M}_{k,n}\times Q$, that the natural map
$\SS\mapor{p}\bar{M}_{k,n}$ is flat and commute with base change; in 
particular 
$p^{-1}(m)=S_m$ and if $E_i\subset \SS$ is the exceptional divisor over 
the section
$p_i$ then $E_i\cap S_m$ is exactly the $(-1)$ curve over the point 
$p_i$.\acapo 
Denoting by $q\colon S\to Q$ the natural projection,  
if $r,s>\sum(a_i+1)$ then the line bundle $\L=q^*\O_Q(r,s)-\sum a_iE_i$ 
is combinatorially ample on every fibre of $p$ and $p_*\L$ is locally free.
\acapo

\Parag{3}{$\zduer$-covers and their deformations}
A quite useful machinery in the explicit construction of components of 
moduli spaces is the technique of Galois deformations of abelian 
covering. This technique has been extensively studied in [Ca1], [Ma2] in 
the case of $\zdue{2}$-covers, in [Ma3] for simple iterated double covers 
and in [Par], [FP] for general nonsingular abelian covers.\acapo
It is possible to prove that the main results of [FP] holds, with minor 
changes, also for flat, locally simple, normal abelian covering. Here, in 
order to avoid long and unnecessary computations we consider only the 
case of $\zduer$-covers, leaving the possible generalizations 
elsewhere.\medskip\acapo 
A) Basic theory of normal flat $\zduer$-covers.\smallskip\acapo

For readers convenience we recall the basic definition and results about 
$\zduer$-covers of algebraic varieties; for the proofs we refer to the 
``standard'' reference [Par].\acapo
It is notationally convenient to consider the group $G=\zduer$ as a vector 
space of dimension $r$ over $\interi/2$; there exists a natural 
isomorphism between the dual vector space $G\dual$ and the group of 
characters $G^*$
$$G\dual\ni\chi\quad\leftrightarrow\quad (-1)^\chi\in G^*.$$
Given any abelian group $\Lambda$, a pair of maps $D\colon G\to \Lambda$, 
$L\colon G\dual\to \Lambda$ satisfies the {\it cover condition} if:\acapo
\item{i)} $D_0=L_0=0$.\acapo
\item{ii)} For every $\chi,\eta\in G\dual$, 
$$L_{\chi}+L_{\eta}=L_{\chi+\eta}+\sum_{\chi(\sigma)=\eta(\sigma)=1}D_\sigma.$$

\acapo
If $i\colon\{0,1\}\to \re$ is the inclusion and $[~]\colon\re\to \interi$ 
is 
the integral part, it is easy to see that a pair 
$(L,D)\in \Lambda^{G\dual}\oplus \Lambda^{G}$ satisfies the cover 
condition if and 
only if $D_0=L_0=0$ and for every $k\ge 0$, $\chi_1,...,\chi_k\in G\dual$ 
holds the equality
$$L_{\chi_1}+\ldots+L_{\chi_k}=L_{\chi_1+...+\chi_k}
+\sum_{\sigma\in G}\Big[{1\over 2}\sum_{j=1}^k i(\chi_j(\sigma))\Big] 
D_\sigma.$$
As an application of this formula we see easily that, 
if $\chi_1,...,\chi_r$ is a basis of $G\dual$, $D\colon 
G\to \Lambda$ is an application such that $D_0=0$ and 
$L_{\chi_1},...,L_{\chi_r}\in \Lambda$ satisfy the equations 
$$2L_{\chi_i}=\sum_{\chi_i(\sigma)=1}D_\sigma,\qquad i=1,...,r$$
then there exists unique an extension $L\colon G\dual\to \Lambda$ such 
that $(L,D)$ satisfies the cover condition.\acapo
We denote by $Z(G,\Lambda)\subset \Lambda^{G\dual}\oplus \Lambda^{G}$ 
the subgroup 
consisting of pairs $(L,D)$ satisfying the cover condition.\acapo
\Lemma{3.1}{} {\it Let $G$ be as above and let
$$0\mapor{}K\mapor{}\tilde{\Lambda}\mapor{\alpha}\Lambda\mapor{}C$$
be an exact sequence of abelian groups with $K$ a $\raz$-vector space and 
$C$ torsion free.\acapo 
Then for every $(L,D)\in Z(G,\Lambda)$ and every 
lifting $\tilde{D}\colon G\to \tilde{\Lambda}$ of $D$ 
there exists unique a lifting 
$(\tilde{L},\tilde{D})\in Z(G,\tilde{\Lambda})$ of $(L,D)$.}\acapo
\Proof It is immediate to prove that if $\alpha(\tilde{\lambda})=2v$ then 
there exist unique $\tilde{v}\in\tilde{\Lambda}$ such that 
$\alpha(\tilde{v})=v$ and 
$2\tilde{v}=\tilde{\lambda}$.\acapo 
Let $\chi_1,...,\chi_r$ be a basis of $G\dual$, by the above remark there 
exists a bijection between the liftings of $L$ satisfying the cover 
conditions and the liftings of $L_{\chi_j}$ satisfying the equations
$$2\tilde{L}_{\chi_i}=\sum_{\chi_i(\sigma)=1}\tilde{D}_\sigma,\qquad 
i=1,...,r.$$
~\finedim
It is very easy to find solutions of the cover condition for every group 
$\Lambda$. 
For every subspace $0\not=H\subset G$ and every $v\in \Lambda$ we can 
define an 
elementary solution 
$(L,D)_{H,v}\in Z(G,\Lambda)$ in the following way (cf. [Par, Ex. 
4.1]):\acapo
i) If $\dim H=1$ we set 
$$L_\chi=\cases{0& if $\chi\in H^{\perp}$,\cr v &if $\chi\not\in 
H^{\perp}$\cr}
\qquad D_\sigma=\cases{2v& if $\sigma\in H-\{0\}$,\cr 0 &otherwise\cr}$$
ii) If $\dim H=d+2\ge 2$
$$L_\chi=\cases{0& if $\chi\in H^{\perp}$,\cr 2^{d}v &if $\chi\not\in 
H^{\perp}$\cr}
\qquad D_\sigma=\cases{v& if $\sigma\in H-\{0\}$,\cr 0 &otherwise\cr}$$
It is not difficult to see that if $\Lambda$ is finitely generated then 
the elements $(L,D)_{H,v}$ generate $Z(G,\Lambda)$; we don't need this 
result.\acapo 
\Lemma{3.2}{} {\it Assume the dimension of $G$ at least 4 and let 
$H\subset G$ a proper subspace.
Let $\alpha\colon \Lambda\to\interi$ be a nonzero  
homomorphism of groups. 
Given any positive integer $N$ and any 
application $D\colon H\to \Lambda$ such that $D_0=0$, there 
exists an extension $(L,D)\in Z(G,\Lambda)$ such that:\acapo 
\item{a)} $\alpha(L_{\chi}-D_\sigma)\ge N$ for every $\chi\in 
G\dual-\{0\}$, $\sigma\in G$.\acapo
\item{b)} $\alpha(D_\sigma)\ge N$ for every $\sigma\in G-H$.\acapo}\acapo
\Proof Let $q$ be a positive integer, $v\in \Lambda$ such that 
$\alpha(v)>0$ and $\eta\in G-H$.\acapo
For every $\tau\in G$ let $(L^\tau, D^\tau)\in Z(G,\Lambda)$ be the 
following elementary solution:\acapo
i) If $\tau\in H$ then
$$L_\chi^\tau=\cases{0& if $\chi\in \eta^{\perp}\cap\tau^{\perp}$,\cr 
D_\tau &otherwise\cr}
\qquad D_\sigma^\tau=\cases{D_\tau& if $\sigma=\tau,\,\eta,\,\tau+\eta$,
\cr 0 &otherwise\cr}$$
ii) If $\tau\not\in H$ then
$$L_\chi^\tau=\cases{0& if $\chi(\tau)=0$,\cr qv &if $\chi(\tau)=1$\cr}
\qquad D_\sigma^\tau=\cases{2qv& if $\sigma=\tau$,
\cr 0 &otherwise\cr}$$
A simple computation shows that for $q>>0$ the sum 
$\ds\sum_{\tau\in G}(L^{\tau}, D^\tau)$ satisfies the required 
conditions.\finedim
Assume now that the group $G$ acts faithfully on a normal irreducible 
complex algebraic variety $X$ and let $\pi\colon X\to Y$ be the 
projection to the quotient. We also make the assumption that the 
quotient map $\pi$ is flat; this  assumption is in general quite strong 
but it is always satisfied  if $Y$ is smooth.\acapo
For every $\sigma\in G-\{0\}$ 
let $Fix(\sigma)\subset X$ be the (possibly empty) closed subscheme of 
fixed points of the involution $\sigma$ and let $R_\sigma\subset 
Fix(\sigma)$ 
be its Weil divisorial part; by convention we set $R_0=\vuoto$. It is 
well known and easy to see (cf. [Ca1]) that $R_\sigma$ is reduced; as the 
group $G$ is abelian every $R_\sigma$ is $G$-stable and if 
$\tau\not=\sigma$ then  $R_\sigma$, $R_\tau$ have no common 
components.\acapo
Let $D_\sigma=\pi(R_\sigma)\subset Y$ with the reduced structure, 
$D_\sigma$ is a Weil divisor which is not Cartier in general.\acapo
The map $\pi$ is flat and then $\pi_*\O_X$ is locally free and there 
exists a character decomposition 
$$\pi_*\O_X=\bigoplus_{\chi\in G\dual}L_{\chi}^{-1}$$
where by definition
$L_{\chi}^{-1}=\{f\in \pi_*\O_X|\,\sigma f=(-1)^{\chi(\sigma)}f,\,   
\forall \sigma\in G\}$, note moreover that $L_0^{-1}=\O_Y$ and, since 
$\pi$ is a $G$-principal bundle over the generic point of $Y$, 
$L_{\chi}^{-1}$ is a locally free $\O_Y$-module of rank 1 for every 
$\chi\in G\dual$.\acapo
For every $x\in X$ let 
$I_x=\{\sigma|\, x\in R_\sigma\}$, 
$I=\cup I_x=\{\sigma|\, R_\sigma\not=\vuoto\}$; 
in the terminology of [Par],  $\pi\colon X\to Y$ is called 
a $(G,I)$-cover.
It is clear that $I_x$ is contained in the stabilizer $Stab_x$ of the 
point $X$; we shall see that the condition $\pi$ flat implies that $I_x$ 
generates $Stab_x$.\acapo
\Def{3.3} In the above notation the cover $\pi\colon X\to Y$ is 
called:\acapo
\item{a)} Totally ramified if $I$ generates $G$.\acapo
\item{b)} Simple if $I$ is a basis of $G$.\acapo
\item{c)} Locally simple if $I_x$ is a set of linearly independent 
vectors for every $x\in X$.\smallskip\acapo
For every pair $\chi,\eta\in G\dual$ let 
$\beta_{\chi,\eta}\colon L_{\chi}^{-1}\otimes L_{\eta}^{-1}\to 
L_{\chi+\eta}^{-1}$ be the multiplication map, it is immediate to observe 
that the $\beta_{\chi,\eta}$'s induce a structure of commutative 
$\O_Y$-algebra on $\oplus L_{\chi}^{-1}$ if and only if 
$$\beta_{\chi,\eta}=\beta_{\eta,\chi},\qquad \beta_{\chi,0}=1,\qquad
\beta_{\chi,\eta}\,\beta_{\chi+\eta,\mu}
=\beta_{\chi,\mu}\,\beta_{\chi+\mu,\eta}\eqno{(3.4)}$$ 
for every $\chi,\eta,\mu\in G\dual$.\acapo
We can interpret $\beta_{\chi,\eta}$ as a section  of the line bundle 
$L_{\chi}+L_{\eta}-L_{\chi+\eta}$ and the computation of [Par] shows that 
$$div(\beta_{\chi,\eta})=\sum_{\chi(\sigma)=\eta(\sigma)=1}D_\sigma$$
Therefore    
over the open set $U\subset Y$ of regular points the divisors $D_\sigma$ 
are Cartier and the pair $(L,D)\in Pic(U)^{G\dual}\oplus Pic(U)^G$ 
satisfies the cover condition. The pair $(L,D)$ is called the {\it 
building 
data} of the cover $\pi\colon X\to Y$.
\acapo
\Cor{3.5}{} {\it Let $\pi\colon X\to Y$ be a normal flat $G$-cover, $Z$ 
an irreducible normal variety and $f\colon Z\to Y$ be a morphism such that 
$f(Z)\not\subset\cup D_\sigma$ and the codimension in $Z$ of 
$f^{-1}(Sing(Y)\cup\pi(Sing(X)))$ is at least 2.\acapo
If the fibred product $X\pr=X\times_YZ$ 
is normal then the divisors 
$f^*(D_\sigma)$ are reduced without common components.}\acapo
\Proof Denote by $\pi\pr\colon X\pr=X\times_YZ\to Z$ the pullback of the 
projection and by 
$f^*\beta_{\chi,\eta}$ the multiplication map of the 
$\O_Z$-algebra $f^*\pi_*\O_X=\pi\pr_*\O_{X\pr}$, then we have 
$$div(f^*\beta_{\chi,\eta})=\sum_{\chi(\sigma)=\eta(\sigma)=1}f^*D_\sigma$$
and the normality of $X\pr$ implies that these divisors are 
reduced.\finedim 
Note that if $\pi$ is locally simple then the divisors $D_\sigma$ are 
Cartier; in fact if $I_x=\{\sigma_1,..,\sigma_s\}$ is a set of linearly 
independent vectors then there exists $\chi_1,...,\chi_s\in G\dual$ such 
that $\chi_i(\sigma_j)=\delta_{ij}$ and therefore by the cover condition 
$2L_{\chi_i}=D_{\sigma_i}$ in some Zariski neighbourhood of 
$y=\pi(x)$.\acapo 
Conversely, given a normal variety $Y$, line bundles $L_{\chi}\in Pic(Y)$ 
and reduced Weil divisors $D_\sigma$ without common components 
such that $(L,D)$ satisfies the cover condition on the Picard group 
of $Y-Sing(Y)$,
there 
exist a normal, flat $G$-cover $\pi\colon X\to Y$ with building data 
exactly $(L,D)$. Moreover if $Y$ is complete then $\pi$ is uniquely 
determined up to $G$-isomorphism.\acapo
In fact let $U\subset Y$ be the open set of smooth points, $\pi\colon 
V\to Y$ be the total space of the vector bundle $\oplus L_{\chi}$, 
$w_{\chi}\in H^0(V,\pi^*L_{\chi})$ the tautological sections and 
$f_\sigma\in H^0(U,D_\sigma)$ a section defining $D_\sigma$.
For  every $\chi,\eta\in G\dual$ let 
$\beta_{\chi,\eta}=\prod_{{\chi(\sigma)=\eta(\sigma)=1}}f_\sigma$; as the 
variety $Y$ is normal, $\beta_{\chi,\eta}$ extends uniquely to a section of
$L_{\chi}+L_{\eta}-L_{\chi+\eta}$. 
An elementary computation shows that the sections $\beta_{\chi,\eta}$ 
satisfy the equations $(3.4)$.\acapo
We then define $X\subset V$ as the 
subvariety defined by the equations 
$$F_{\chi,\eta}:=w_\chi w_\eta-w_{\chi+\eta}\beta_{\chi,\eta}=0,\quad 
\forall \chi,\eta\in G\dual.$$
For the unicity of $\pi$ when $Y$ is complete 
follows from the following argument; another proof will be made 
later on this section.\acapo
Let $f\pr_\sigma$ be another set of equation of the divisors $D_\sigma$ 
and let $\beta_{\chi,\eta}\pr$ be the corresponding multiplication maps 
and $X\pr\subset V$ defined by the equations 
$w_\chi w_\eta-w_{\chi+\eta}\beta_{\chi,\eta}\pr=0$.
For every pair $\chi,\eta\in G\dual$ there exists an invertible function 
$\alpha_{\chi,\eta}\in H^0(\O_Y^*)$ such that 
$\beta_{\chi,\eta}=\alpha_{\chi,\eta}\beta_{\chi,\eta}\pr$.\acapo
Assume now that for every $\chi\in G\dual$ there exists a square root of 
$\alpha_{\chi,\chi}$ in $H^0(\O_Y^*)$, then there exist 
invertible functions $g_\chi\in H^0(\O_Y^*)$ such that 
$$g_\chi g_\eta=g_{\chi+\eta}\alpha_{\chi,\eta}\eqno{(3.6)}$$ 
for every $\chi,\eta\in 
G\dual$; the isomorphism of the vector bundle $V$ defined by $w_\chi\to 
g_\chi w_\chi$ will induces an isomorphism $X\simeq X\pr$.\acapo
We can choose the function $g_\chi$ in the following way; let 
$\chi_1,...,\chi_r$ be a basis of $G\dual$ and for every $\chi=\sum 
a_i\chi_i$, $a_i=0,1$, we define the length $l(\chi)=\sum a_i$. Let 
$g_{\chi_i}$ be a fixed square root of $\alpha_{\chi_i,\chi_i}$;
let $\chi=\eta+\chi_i$ be a decomposition such that $l(\chi)=l(\eta)+1$, 
by induction on $l$ we can suppose $g_\eta$ is defined, then we set 
$g_\chi=g_\eta g_{\chi_i}\alpha_{\eta,\chi_i}^{-1}$: the functions 
$g_\chi$ satisfy the cocycle condition $(3.6)$.\acapo
\Prop{3.7}{} {\it Let $\pi\colon X\to Y$ be a normal flat $G$-cover with 
$Y$ 
smooth. Then $X$ is smooth if and only if $\pi$ is locally simple, 
the divisors $D_\sigma$ are 
smooth and intersects transversally.}\acapo
\Proof [Par,3.1].\finedim
In the case $Y$ smooth surface there exists simple formulas for the 
numerical invariants $\chi(\O_X)$, $K^2_X$ in terms of the building data 
$(L,D)$. In fact by Hurwitz formula we have $K_X=\pi^*K_Y+\sum R_\sigma$, 
$\pi^*D_\sigma=2R_\sigma$ and if $g=|G|=2^r$ then for every divisors $A,B$ 
on $Y$ we have 
$$\pi^*A\per\pi^*B=gA\per B,\quad R_\sigma\per \pi^*A={1\over 
2}gD_\sigma\per A,\quad R_\sigma\per R_\tau={1\over 4}g D_\sigma\per 
D_\tau$$
Since $\pi$ is finite $H^i(\O_X)=\oplus H^i(Y, L_{\chi}^{-1})$ and 
therefore we get 
$$K_X^2=g(K_Y+{1\over 2}\sum_\sigma D_\sigma)^2,\qquad 
\chi(\O_X)=g\chi(\O_Y)+{1\over 2}\sum_{\chi}L_{\chi}\per(L_{\chi}+K_Y)$$
A simple calculation shows that if the building data are sufficiently 
ample then the invariants $K^2,\chi$ spread in the region $4\chi\le 
K^2\le 8\chi$.\acapo 
\Rem{} When we say that the building data $(L,D)$ is sufficiently ample 
we shall mean that there exists a basis $\sigma_1,...,\sigma_r$ of $G$ 
such that the line bundles $L_{\chi}$, $\chi\not=0$ and 
$\O_Y(D_{\sigma_i})$, $i=1,..,r$, are sufficiently ample.\acapo
By the computation of [Ma3] it follows that every real number $\beta\in 
[4,8]$ can be approximated by the ratios $\ds{K^2\over\chi}$ of simple 
$\zduer$-covers.
\medskip\acapo
B) The geometry of locally simple covers.\smallskip\acapo
Let $\pi\colon X\to Y$ be a normal flat simple (G,I)-cover with 
$I=\{\sigma_1,...,\sigma_r\}$ basis of $G=\zduer$ and $Y$ complete 
irreducible. 
Let 
$\chi_1,...,\chi_r\in G\dual$ be the dual basis of 
$\sigma_1,...,\sigma_r$ and let $f_i$ be an equation defining 
$D_{\sigma_i}$. In the above notation $X$ is defined, up to isomorphism, 
in $V$ by the equations 
$$\cases{w_{\chi_i}^2=f_i\quad & for $i=1,...,r$\cr \quad &\cr
w_{\sum a_i\chi_i}=\prod w_{\chi_i}^{a_i}\quad & for $a_i=0,1$\cr}$$
An immediate simplification of this equations gives $X$ as a complete 
intersection in the total space of $L_{\chi_1}\oplus....\oplus L_{\chi_r}$ 
defined by the $r$ equations $w_{\chi_i}^2=f_i$.
Note moreover that the divisor $R_{\sigma_i}$ is defined by the equation 
$w_{\chi_i}=0$ and then it is a Cartier divisors in $X$.
From  this description follows immediately the projection formula 
$$\pi_*\O_{R_\sigma}=\Somdir{\chi(\sigma)=0}{}\O_{D_\sigma}(-L_\chi).$$ 
Note that our definition of simple cover is 
compatible with the, apparently different, 
notions of simple covers given in [Ma3], [Ca3], [Par].\acapo
If $\pi$ is locally simple then for every $y\in Y$ there exists a  
neighbourhood $U\subset Y$ such that $X\times_Y U$ is isomorphic to the 
cover 
defined by equations $w_i^2=f_i$, $i=1,...,r$. In particular 
a normal flat $G$-cover is 
unramified at $x$ if and only if $I_x=\vuoto$; more generally if $H\subset 
G$ 
is the subspace generated by $I_x$ then $X/H\to Y$ if still flat and 
ramified in codimension $\ge 2$ at the point $x$, by the above  remark 
$X/H\to Y$ is unramified at $x$ and then $H=Stab_x$.\acapo
A more important consequence is that if $\pi\colon X\to Y$ is locally 
simple then  
$X$ is a local complete intersection in $V$ and the divisors $R_\sigma$ 
are Cartier.\acapo
Let $\pi\colon X\to Y$ be a locally simple $G$-cover, we have seen that 
for every $\sigma$, $R_\sigma$ is a $G$-stable Cartier divisors, in 
particular the $G$-action on $X$ can be naturally lifted, in the sense of 
[Mu2, p.110-111], to $G$-actions on the coherent sheaves
$\O_X(-R_\sigma)$, $\O_{R_\sigma}(-R_\sigma)$.\acapo
\Theor{3.8}{} {\it For every locally simple normal flat $\zduer$-cover 
$\pi\colon X\to Y$ there exists a $G$-equivariant exact sequence of sheaves
$$0\mapor{}\pi^*\Omega^1_Y\mapor{}\Omega^1_X\mapor{}
\oplus_{\sigma}\O_{R_\sigma}(-R_\sigma)\mapor{}0.$$}\acapo
\Proof Assume first $\pi$ simple, then there exist line bundles 
$L_1,...,L_r$ over $Y$ such that $X$ is the complete intersection 
subvariety of $W=L_1\oplus...\oplus L_r$ defined by the equations
$F_i=w_i^2-f_i=0$, where $D_i=\{f_i=0\}\subset Y$ are the branching 
divisors of 
$\pi$ and $w_i\in H^0(W,\pi^*L_i)$ are the tautological sections.\acapo
The equations $F_i\in H^0(W,\pi^*D_i)$ 
are $G$-invariant, in particular the conormal bundle 
of $X$ in $W$ is $G$-isomorphic to the pull back 
$\pi^*(\oplus\O_Y(-D_i))=\oplus\O_X(-2R_i)$. 
On the other hand there exists a natural 
$G$-isomorphism $\Omega^1_{W/Y}=\oplus \pi^*L_i^{-1}$. Tensoring the 
exact sequence of sheaves over $W$  
$$0\mapor{}\pi^*\Omega^1_{Y}\mapor{} 
\Omega^1_W\mapor{}\Omega^1_{W/Y}\mapor{}0$$ 
by $\O_X$ we get 
$$Tor_1^W(\O_X,\Omega^1_{W/Y})=0\mapor{}\pi^*\Omega^1_{Y}\mapor{} 
\Omega^1_W\otimes\O_X\mapor{}\Omega^1_{W/Y}\otimes\O_X\mapor{}0.$$
We have moreover the exact sequence 
$$0\mapor{}\oplus_i\O_X(-2R_i)\mapor{d} 
\Omega^1_{W/Y}\otimes\O_X\mapor{}\Omega^1_{X/Y}\mapor{}0,$$ 
where $d(F_i)=d(w_i^2-f_i)=w_idw_i\in\O_X(-R_i)=\pi^*L_i^{-1}$. 
This gives a natural 
$G$-isomorphism  
$$\Omega^1_{X/Y}=\bigoplus_i{\O_X(-R_i)\over 
w_i\O_X(-2R_i)}=\bigoplus_i\O_{R_i}(-R_i).$$ 
The injectivity of $\pi^*\Omega_Y^1\to\Omega^1_X$ follows by considering 
the snake exact sequence of the commutative diagram
$$\matrix{0&\mapor{}&\oplus\O_X(-\pi^*D_i)&\mapor{}&\Omega_W^1\otimes\O_X&
\mapor{}&\Omega^1_X&\mapor{}&0\cr
&&\|&&\mapver{}&&\mapver{}&&\cr
0&\mapor{}&\oplus\O_X(-\pi^*D_i)&\mapor{}&\Omega_{W/Y}^1\otimes\O_X&
\mapor{}&\Omega^1_{X/Y}&\mapor{}&0\cr}$$
This proves the theorem for simple covers.\acapo
In the general case, for every $\sigma\in G$ let $Z_\sigma=X/\sigma$ and 
$\pi_\sigma\colon X\to Z_\sigma$  be the projection to the quotient. As 
$R_\sigma$ is Cartier, by [Ma3] 3.1, the map $\pi_\sigma$ is flat in a 
neighbourhood $U$ of $R_\sigma$ and by the  computation in the simple 
case there exists a natural isomorphism 
$\Omega^1_{U/Z_\sigma}=\O_{R_\sigma}(-R_\sigma)$ and morphisms 
$$\Omega_{X/Y}^1\to\Omega^1_{X/Z_\sigma}\to \O_{R_\sigma}(-R_\sigma)$$
Taking the direct sum of $\O_{R_\sigma}(-R_\sigma)$ over all $\sigma\in 
G$ we get a morphism $\Omega^1_{X/Y}\to \oplus\O_{R_\sigma}(-R_\sigma)$ 
and the computation made in the simple case shows that it is a 
$G$-isomorphism.\finedim 
\Lemma{3.9}{} {\it Let $\pi\colon X\to Y$ be a normal flat 
locally simple $G$-cover, then 
for every $\sigma\in G$
$$\Ext^i_X(\O_{R_\sigma}(-R_\sigma),\O_X)=
\bigoplus_{\chi(\sigma)=0}H^{i-1}(\O_{D_\sigma}(D_\sigma-L_\chi))$$}
\Proof By applying the functor $\HOM_{X}(-,\O_X)$ to the exact sequence 
$$0\mapor{}\O_X(-2R_\sigma)\mapor{}\O_X(-R_\sigma)
\mapor{}\O_{R_\sigma}(-R_\sigma)\mapor{}0$$ 
we get
$$\EXT^i_X(\O_{R_\sigma}(-R_\sigma),\O_X)=
\cases{\O_{R_\sigma}(2R_\sigma)\qquad &$i=1$\cr
0\qquad &$i\not=1$\cr}$$
and by the $\Ext$ spectral sequence 
$\Ext^i_X(\O_{R_\sigma}(-R_\sigma),\O_X)=
H^{i-1}(\pi_*\O_{R_\sigma}(2R_\sigma))$.\acapo
A local computation shows that $2R_\sigma=\pi^*D_\sigma$, 
$\pi_*\O_{R_\sigma}=\oplus_{\chi(\sigma)=0}\O_{D_\sigma}(-L_\chi)$ and 
then 
$$\pi_*\O_{R_\sigma}=\bigoplus_{\chi(\sigma)=0}\O_{D_\sigma}(-L_\chi).$$
~\finedim
\Cor{3.10}{} {\it If $\pi\colon X\to Y$ is normal, flat, locally simple 
$G$-cover such that 
$H^1(Y,L_\chi^{-1})=0$, $\Ext^1_Y(\Omega^1_Y,L_\chi^{-1})=0$, 
$H^0(Y,D_\sigma-L_\chi)=0$ for every $\chi\not=0$ and every $\sigma\in 
\chi^{\perp}$.\acapo 
Then $G$ acts trivially on $\Ext^1_X(\Omega_X^1,\O_X)$.}\acapo
\Proof Let $\chi\not=0$ be a fixed character, taking the 
$\chi$-equivariant part of the long $\Ext_X(-,\O_X)$ sequence associated 
to 
$$0\mapor{}\pi^*\Omega^1_Y\mapor{}\Omega^1_X\mapor{}
\oplus_{\sigma}\O_{R_\sigma}(-R_\sigma)\mapor{}0$$
we get 
$$\bigoplus_{\sigma\in\chi^{\perp}}H^0(\O_{D_\sigma}(D_\sigma-L_\chi))
\mapor{}\Ext^1_X(\Omega_X^1,\O_X)^\chi\mapor{}
\Ext^1_X(\pi^*\Omega_Y^1,\O_X)^\chi$$
The left side lies in the middle of $H^0(Y,D_\sigma-L_\chi)$ and 
$H^1(Y,-L_\chi)$ and then it is equal to 0. As $\pi$ is finite flat there 
exists a natural isomorphism (cf. [Ma3] 2.1) 
$\Ext^1_X(\pi^*\Omega_Y^1,\O_X)=\Ext^1_Y(\Omega_Y^1,\pi_*\O_X)$ and 
therefore 
$\Ext^1_X(\pi^*\Omega_Y^1,\O_X)^\chi=\Ext^1_Y(\Omega_Y^1,L_{\chi}^{-1})=0$.
\finedim 
Note that if $Y$ is Gorenstein of dimension $n\ge 2$ then by Serre duality 
$$\Ext^1_Y(\Omega_Y^1,L_{\chi}^{-1})=
\Ext^1_Y(\Omega_Y^1\otimes L_\chi\otimes K_Y,K_Y)=
H^{n-1}(\Omega_Y^1\otimes L_\chi\otimes K_Y)\dual$$
and this space is $0$ whenever $L_\chi$ is sufficiently 
ample.\bigskip\acapo
C) Galois deformations of $\zduer$-covers.\medskip\acapo
Let $\pi\colon X\to Y$ be a normal flat $\zduer$-cover with $Y$ complete 
irreducible variety 
and building data 
$L_\chi$, $D_\sigma$. Here we want to study the deformations of $\pi$ 
obtained by ``moving'' the branching divisors.
We assume for simplicity that 
$Y$ is smooth and $H^0(Y,\theta_Y)=H^0(X,\theta_X)=0$; although these 
assumption are not strictly necessary for the main results of this 
section, 
they are 
sufficient for our application and allow simpler proofs. We also note 
that in the case $X$ smooth our results are  contained in 
[FP].\acapo 
Let $Art$ be the (small) category of local Artinian $\com$-algebras and 
denote by 
$$Def_X,Def_Y,Def_\pi\colon Art\to Set$$
the functors of deformations of $X$, $Y$, $\pi$ respectively. 
All of them satisfy Schlessinger conditions H1, H2, H3 of [Sch] and 
linearity condition L of [FaMa]; therefore they admit a hull and a good 
obstruction theory.\acapo
Since $X,Y$ 
are normal we have 
$$T^iDef_X=\Ext^i_X(\Omega^1_X,\O_X),\qquad 
T^iDef_Y=\Ext^i_Y(\Omega^1_Y,\O_Y),\qquad i=1,2$$
where as usual we denote by $T^1$ the tangent space and by $T^2$ the 
obstruction space arising from the cotangent complex [Fle].\acapo
Let $Def_{(Y,D_\sigma)}\colon Art\to Set$ be the functor of deformations 
of the closed inclusions $D_\sigma\to Y$; more precisely for $A\in Art$, 
$Def_{(Y,D_\sigma)}(A)$ is the set of isomorphism classes of:
\item{1)} a deformation of $Y$, $Y_A\to Spec(A)$.\acapo
\item{2)} for every $\sigma\in G$, a 
closed $A$-flat embedding $D_{A,\sigma}\subset Y_A$ extending 
$D_\sigma$ (the subvarieties $D_{A,\sigma}$ will be automatically Cartier 
divisors).\acapo 
Note that $Def_{(Y,D_\sigma)}$ is prorepresented by the fibred product of 
the relative Hilbert scheme of the Kuranishi family of $Y$.\acapo
\Lemma{3.11}{} {\it In the above set up let $(Y_A,D_{A,\sigma})\in 
Def_{(Y,D_\sigma)}(A)$; then for every $\chi$ there exists a unique 
extension $L_{A,\chi}\in Pic(Y_A)$ of $L_\chi$ such that 
$(L_{A,\chi},D_{A,\sigma})$ satisfies the cover condition.}\acapo
\Proof Let $0\to \com\to A\to B\to 0$ be a small extension in $Art$, then 
we have an exact sequence of sheaves
$$0\mapor{}\O_Y\mapor{}\O^*_{Y_A}\mapor{}\O_{Y_B}^*\mapor{}.0$$
Since $Y$ is complete, reduced and irreducible 
the pullback map $B^*=(B-m_B)\to H^0(\O_{Y_B}^*)$ is surjective 
and then we have an exact sequence of abelian groups 
$$0\mapor{}H^1(\O_Y)\mapor{}Pic(Y_A)\mapor{}Pic(Y_B)\mapor{}H^2(\O_Y).$$
The conclusion follows by lemma 3.1 and induction on the length of 
$A$.\finedim 
\Prop{3.12}{} {\it There exists a natural transformation of functors
$$Def_{(Y,D_\sigma)}\to Def_\pi$$
commuting with the projections $Def_\pi\to Def_Y$, 
$Def_{(Y,D_\sigma)}\to Def_Y$.\acapo}
\Proof 
Let $f_\sigma\in H^0(Y,D_\sigma)$ be equations of the divisors $D_\sigma$, 
then $X$ is $G$-isomorphic to the subvariety of $V=\oplus L_\chi$ defined 
by 
$$F_{\chi,\eta}=w_\chi 
w_\eta-w_{\chi+\eta}\prod_{\chi(\sigma)=\eta(\sigma)=1}f_\sigma=0,\qquad
\qquad\chi,\eta\in G\dual.$$ 
Consider now $A\in Art$ and $(Y_A,D_{A,\sigma})\in 
Def_{(Y,D_\sigma)}(A)$. By 3.11 there exist unique line bundles 
$L_{A,\chi}\in Pic(Y_A)$ such that $(L_{A,\chi},D_{A,\sigma})$ satisfies 
the cover condition.\acapo
Denoting $\pi_A\colon V_A=\oplus L_{A,\chi}\to Y_A$ 
we define 
$X_A\subset V_A$ by the equations 
$$F_{A,\chi,\eta}=w_\chi w_\eta-w_{\chi+\eta}\prod_{\chi(\sigma)
=\eta(\sigma)=1}f_{A,\sigma}=0\qquad\qquad\chi,\eta\in G\dual.$$ 
where $f_{A,\sigma}\in H^0(Y_A,D_{A,\sigma})$ is an equation of 
$D_{A,\sigma}$ extending $f_\sigma$. We also define $\pi_A\colon X_A\to 
Y_A$ as the natural projection.\acapo
We first note that $\pi_A$ is flat and then $X_A$ is flat over $A$. In 
fact by construction $\pi_{A*}\O_{X_A}=\oplus L_{A,\chi}^{-1}$ which is 
locally free. It is clear that our construction commutes with morphisms in 
$Art$ and then remains only to prove that the isomorphism deformation 
class of 
$\pi_A$ does not depends on the choice of $f_{A,\sigma}$.
More precisely we also need to prove the independence from the choice of 
the isomorphisms of line bundles $L_{A,\chi}\otimes L_{A,\eta}\otimes 
L_{A,\chi+\eta}^{-1}\simeq 
\O(\sum_{\chi(\sigma)=\eta(\sigma)=1}D_{A,\sigma})$, but the effect over 
$\pi_A\colon X_A\to Y_A$ of a 
change of isomorphism is the same of a change of the equations 
$f_{A,\sigma}$.\acapo 
Let $f_{A,\sigma}\pr$ be another set of equations, as $Y$ is complete 
reduced and irreducible, for every $\sigma$ there exists an invertible 
element $a_\sigma\in A$ such that $f_{A,\sigma}\pr=a_\sigma 
f_{A,\sigma}$. Let $\pi_A\pr\colon X_A\pr\to Y_A$ be the cover defined by 
$f_{A,\sigma}\pr$.\acapo 
The maps $\pi_A,\pi_A\pr$ are uniquely determined by the isomorphism 
classes of the $\O_{Y_A}$-algebras 
$$\A=\oplus L_{A,\chi}^{-1},\qquad 
\A\pr=\oplus L_{A,\chi}^{-1}$$ 
with the 
multiplication maps induced respectively by the sections
$$\beta_{A,\chi,\eta}=
\prod_{\chi(\sigma)=\eta(\sigma)=1}f_{A,\sigma},\qquad
\beta_{A,\chi,\eta}\pr=
\prod_{\chi(\sigma)=\eta(\sigma)=1}f_{A,\sigma}\pr$$
The sections $\beta,\beta\pr$ differ by multiplication by an invertible 
element of $A$.\acapo
Every $\O_{Y_A}$-isomorphism from $\A$ to $\A\pr$ is determined by taking 
for every $\chi\in G\dual$ an invertible element $g_\chi\in 
\Hom(L_{A,\chi},L_{A,\chi})=A$ such that $g_0=1$ and 
$$\beta_{A,\chi,\eta}\pr={g_\chi g_\eta\over 
g_{\chi+\eta}}\beta_{A,\chi,\eta},\qquad\forall \chi,\eta\in G\dual$$
Given a section $s\colon Spec(A)\to Y_A$
it is clearly sufficient to check these conditions over the image of $s$, 
or equivalently to prove that $s^*\A$ is isomorphic to 
$s^*\A\pr$. This last condition is trivially satisfied if the image of 
$s$ does not intersect the branching divisors $D_{A,\sigma}$, in fact 
every principal $G$-bundle over a fat point is trivial.\finedim 
Let $\eta\colon Def_{(Y,D_\sigma)}\to Def_X$ be the composition of the 
morphisms defined in 3.8 and the natural projection $Def_\pi\to Def_X$. 
Let moreover $Def_X^G\subset Def_X$ be the subfunctor of deformations of 
$X$ for which there exists an extension of the $G$-action. Using the fact 
that $H^0(X,\theta_X)=0$ it is easy to see that $Def_X^G$ is a functor 
with good deformation theory, i.e. satisfies the Schlessinger conditions, 
and it is prorepresented by the $G$-invariant subgerm of the Kuranishi 
base space of $X$. By Cartan's lemma $Def_X^G=Def_X$ if and only if $G$ 
acts trivially on $T^1_X=\Ext^1_X(\Omega^1_X,\O_X)$.\acapo
\Theor{3.13}{} {\it In the above notation, if $\pi\colon X\to Y$ is 
locally simple then the natural transformation $\eta$ induces an 
isomorphism 
$Def_{(Y,D_\sigma)}\simeq Def_X^G$.}\acapo
\Proof We prove the theorem by constructing explicitly an inverse natural 
transformation $\theta\colon Def_X^G\to Def_{(Y,D_\sigma)}$. The first 
step is to construct a natural transformation $\mu\colon Def_X^G\to 
Def_Y$ by ``taking the quotients''.\acapo
For every $A\in Art$ it is useful to consider every 
deformation  $X_A\in Def_X^G(A)$ as a sheaf $\O_A$ of flat $A$-algebras 
over $X$ together with a lifting of the $G$-action. As $G$ is abelian 
there exists a character decomposition 
$$\pi_*\O_A=\bigoplus_{\chi\in G\dual}(\pi_*\O_A)^\chi$$
Since $\pi$ is a finite affine morphism, the sheaf $\pi_*\O_A$ is still 
flat over $A$ and then every direct summand $(\pi_*\O_A)^\chi$ is 
$A$-flat. In particular $Y_A\colon=(Y,(\pi_*\O_A)^0)\to Spec(A)$ is a 
flat morphism.\acapo 
If $A\to B$ is a morphism in $Art$ and $\O_B=\O_A\otimes_AB$ we have 
$\pi_*\O_B=(\pi_*\O_A)\otimes_AB$;  moreover
for every character $\chi$ we have 
$(\pi_*\O_A)^\chi\otimes_AB\subset \pi_*(\O_B)^\chi$. Putting together 
the above relations we get 
that $(\pi_*\O_A)^\chi\otimes_AB=\pi_*(\O_B)^\chi$; in particular $Y_A$ 
is a deformation of $Y$ and the induced application 
$\mu\colon Def_X^G\to Def_Y$ commute with base change, i.e. it is a 
morphism of functors.\acapo
Note also that $(\pi_*\O_A)^\chi$ is a deformation of $L_\chi^{-1}$ and 
therefore it is a line bundle over $Y_A$.\acapo
Let $x\in X$ be a point, $y=\pi(x)$ and $\sigma_1,...,\sigma_r\in G$ a 
basis extending $I_x$; denote by $\chi_1,...,\chi_r$ the dual basis. Over 
a sufficiently small affine open neighbourhood $U$ of $y$ the cover 
$\pi\colon X\to Y$ is isomorphic to $\pi_1\colon X_U\to U$ 
where $X_U\subset U\times \com^r_w$ is defined by equations 
$$w_i^2=f_{\sigma_i},\qquad i=1,...,r$$
and $\pi_1$ is the projection on the first factor. The $G$-action over 
$X_U$ is induced by the $G$-action on the affine variety 
$U\times \com^r$ defined by 
$\sigma(w_i)=(-1)^{\chi_i(\sigma)}w_i$.\acapo 
Let $(X,\O_A)\in Def_X^G(A)$ and let $U_A=(U,(\pi_*\O_A)^0)$, it is well 
known [Ar] that $(X_U,\O_A)\to U_A$ 
can be described as an embedded deformation inside $U_A\times\com^r$, 
moreover as $G$ is finite we can choose the embedded deformation 
$G$-stable.\acapo 
$X_U$ is a complete intersection in $U\times\com^r$ and then every 
$G$-stable deformation over $A$ is given by  the equations
$$w_i^2=f_{A,\sigma_i},\qquad i=1,...,r$$
with $f_{A,\sigma_i}$ a lifting of $f_{\sigma_i}$ over $U_A$.\acapo
The regular function $f_{A,\sigma}$ defines a Cartier divisor 
$D_{A,\sigma}$ over $U_A$ which is $A$-flat. It is clear that these local 
constructions of $D_{A,\sigma}$ patch together over $Y_A$ and commute with 
base change. \acapo
We have therefore constructed a natural transformation 
$\theta\colon Def_X^G\to Def_{(Y,D_\sigma)}$. By construction it is 
immediate to observe that $\theta\eta$ is the identity on 
$Def_{(Y,D_\sigma)}$.\acapo 
The same proof of  3.12 shows that the morphism $\theta$ is injective  
and then  $\eta\theta$ is the identity on $Def_X^G$.\finedim
The above proof has the advantage of being quite elementary but it seems 
very hard to extend it to natural deformations ([FP] def. 3.2). 
For this reason we sketch here a shorter but more technical proof of 3.13
which involves the notion of relative obstructions and the standard 
criterion of smoothness [FaMa].\acapo
Let $\nu\colon Def_{(Y,D_\sigma)}\to Def_Y$ be the natural projection, it 
is easy to compute the tangent and obstruction space of $\nu$; they are
$$T^i\nu=\bigoplus_\sigma H^{i-1}(\O_{D_\sigma}(D_\sigma)),\qquad i=1,2.$$
This gives an exact sequence 
$$0\mapor{}\bigoplus_\sigma 
H^{0}(\O_{D_\sigma}(D_\sigma))\mapor{}T^1Def_{(Y,D_\sigma)}
\mapor{}\Ext^1_Y(\Omega^1_Y,\O_Y)\mapor{}$$  
$$\mapor{}\bigoplus_\sigma 
H^{1}(\O_{D_\sigma}(D_\sigma))\mapor{}T^2Def_{(Y,D_\sigma)}
\mapor{}\Ext^2_Y(\Omega^1_Y,\O_Y)$$  
On the other hand, taking the long $\Ext_X(-,\O_X)$ sequence of 
$$0\mapor{}\pi^*\Omega^1_Y\mapor{}\Omega^1_X\mapor{}
\oplus_\sigma\O_{R_\sigma}(-R_\sigma)\mapor{}0$$ 
we get 
$$\bigoplus_\chi \Ext^0_Y(\Omega^1_Y,L_\chi^{-1}) 
\mapor{}\bigoplus_\sigma \bigoplus_{\chi(\sigma)=0}
H^{0}(\O_{D_\sigma}(D_\sigma-L_\chi))\mapor{}T^1Def_X
\mapor{}\bigoplus_\chi\Ext^1_Y(\Omega^1_Y,L_\chi^{-1})\mapor{}$$  
$$\mapor{}\bigoplus_\sigma \bigoplus_{\chi(\sigma)=0}
H^{1}(\O_{D_\sigma}(D_\sigma-L_\chi))\mapor{}T^2Def_X
\mapor{}\bigoplus_\chi\Ext^2_Y(\Omega^1_Y,L_\chi^{-1})$$  
If $\Ext^1_Y(\Omega^1_Y,L_\chi^{-1})=0$, 
$H^{0}(\O_{D_\sigma}(D_\sigma-L_\chi))=0$ 
for every $\chi\not=0$ and $\sigma\in \chi^{\perp}$, then a simple 
diagram chasing involving the above exact sequences 
shows that the natural maps 
$$T^iDef_{(Y,D_\sigma)}\to T^iDef_X$$
are surjective for $i=1$ and injective for $i=2$. By the standard 
smoothness criterion the morphism $\eta$ is smooth. 
If in addition $H^0(\theta_Y)=\Ext^0_Y(\Omega^1_Y,\O_Y)=0$ then
also $H^0(\theta_X)=0$ and $Def_X$ is prorepresentable. By diagram 
chasing  
$\eta$ induces an isomorphism on tangent spaces and then 
it is an isomorphism.\acapo
For future reference we collect the main results in the following corollary
\Cor{3.14}{} {\it Let $\pi\colon X\to Y$ be a normal flat locally simple 
$\zduer$-cover with building data $(L_\chi,D_\sigma)$ such that:\acapo
\item{i)} $Y$ is smooth of dimension $\ge 2$, $H^0(\theta_Y)=0$.\acapo 
\item{ii)} $\Ext^1_Y(\Omega^1_Y,L_\chi^{-1})=H^1(L_\chi^{-1})=0$ for 
every $\chi\not=0$.\acapo 
\item{iii)} $H^0(Y, D_\sigma-L_\chi)=0$ for every $\chi\not=0$, 
$\sigma\in \chi^{\perp}$.\acapo 
Then $H^0(\theta_X)=0$ and the natural morphism  
$\eta\colon Def_{(Y,D_\sigma)}\to Def_X$ is an isomorphism of 
prorepresentable functors.}\acapo
The functor $Def_{(Y,D_\sigma)}$ corresponds essentially to the functor 
$Dgal_X$ of Galois deformations introduced in [FP]; we shall say that 
{\it the family of Galois deformations is complete} if $\eta$ is smooth.\acapo 
Let  $H\subset G$ be a proper 
subspace and let $D_\sigma$, $\sigma\in H-\{0\}$, be  reduced effective
divisors in the smooth variety $Y$ without common components.
The same argument of lemma 3.2, with $\Lambda=Pic(Y)$ and $v$ ample, 
shows that if the dimension of $G$ is at least 4 there exists 
sufficiently ample and general divisors 
$D_\sigma$, $\sigma\in G-H$, and sufficiently ample line bundles
$L_\chi$ such that the pair $(L,D)$ satisfies the cover condition. 
Moreover if $\dim Y\ge 2$ and $H^0(\theta_Y)=0$ we can choose $(L,D)$ 
which satisfies the hypothesis of corollary 3.14. In particular if $X\to 
Y$ is the $G$-cover with building data $(L,D)$ then the family of Galois 
deformations of $X$ is complete.\acapo
\Ex{3.15}{} The corollary 3.14 provides a techniques for constructing 
irreducible components of the moduli space of surfaces of general type 
whose members have singular canonical model.\acapo
Consider for example the Pascal configuration in $\pro^2$, i.e. a smooth 
conic $C$ and an inscribed hexagon $E=L_1+L_2+...+L_6$, $L_i$ line. This 
configuration is uniquely determined by $C$ and the six vertices 
$p_i=C\cap L_i\cap L_{i+1}$, $(L_7=L_1)$, $i=1,...,6$. Assume 
$p_1,...,p_6\in C$ in general position and let $q_1=L_1\cap L_4$, 
$q_2=L_2\cap L_5$, $q_3=L_3\cap L_6$. Pascal theorem asserts that the 
$q_i$ are collinear end then the generic quadric passing through $q_i$ 
has an ordinary double point (node).\acapo
Let $Y\mapor{\nu}\pro^2$ be the blow up of $p_1,...,p_6,q_1,q_2,q_3$. 
By Kodaira stability theorem  
the 9 exceptional curves of $\nu$ are stable under deformations, in 
particular it is possible to choose
a very ample line bundle $\O_Y(1)\in Pic(Y)$ which extends to every 
deformation of $Y$.\acapo
Let $\sigma_1,...,\sigma_8\in \zdue{8}$ be a basis and let 
$D\colon H=\zdue{8}\to Div(Y)$ defined in the following way:\acapo
\item{} $D_{\sigma_i}$= strict transform of $L_i$, $i=1,...,6$.\acapo
\item{} $D_{\sigma_7}=$ strict transform of $C$.\acapo
\item{} $D_{\sigma_8}=$ strict transform of a generic quadric passing 
through $q_1,q_2,q_3$.\acapo
\item{} $D_\sigma=0$ if $\sigma\not=\sigma_i$.\acapo
It is immediate to see that every deformation of the data $(Y,D_\sigma)$ 
preserves the Pascal configuration; in particular if $(Y_A, 
D_{A,\sigma})\in Def_{(Y,D_\sigma)}(A)$ then the flat family of curves 
$D_{A,\sigma_8}\to Spec(A)$ is locally trivial.\acapo
Let now $D\colon G=\zdue{9}\to Div(Y)$, $L\colon G\to Pic(Y)$ 
an extension of $D\colon H\to Div(Y)$ such that 
for every $\sigma\in G-H$ the divisor 
$D_\sigma$ is a generic element of the linear system $|\O_Y(l)|$, 
$l>>0$; 
the pair $(L,D)$ satisfies the cover conditions and the hypothesis of 
Corollary 3.14.\acapo
Let $\pi\colon X\to Y$ be the $G$-cover with building data $(L,D)$, it is 
easy to see that $X$ is a normal surface with exactly $2^8$ 
simple nodes as singularities and ample canonical bundle. 
By 3.14 every deformation of $X$ is Galois and then contains at least 
$2^8$ nodes.\acapo
Examples of irreducible components of the moduli space of surfaces of 
general type whose general member is a regular surface with singular 
canonical model were constructed in [Ca7] by using a different method.
\medskip\acapo
The use of Galois deformations allow to give a simple proof of the 
following 
\Prop{3.16}{} {\it Let $\pi\colon X\to Y$ be a normal flat $\zduer$-cover 
with $Y$ smooth; then $X$ is $\raz$-Gorenstein of index $\le 2$. Moreover 
the following three conditions are equivalent:\acapo
\item{a)} $\pi$ is locally simple.\acapo
\item{b)} $X$ is locally complete intersection.\acapo
\item{c)} $X$ is Gorenstein.\acapo}
\Proof 
As $Y$ is smooth and $\pi$ is flat, 
by a well known result in commutative algebra [Mat,\S 23], the 
variety $X$ is Cohen-Macaulay (resp.: Gorenstein) 
if and only if the fibres of $\pi$ are 
Cohen-Macaulay (resp.: Gorenstein).
The fibres, and hence $X$, are Cohen-Macaulay 
because they have dimension 0. By Hurwitz formula 
$2K_X=f^*K_Y+\sum D_\sigma$; this proves that $X$ is $\raz$-Gorenstein of 
index $\le 2$.\acapo 
We have already proved that $a)\solose b)$, while $b)\solose c)$ 
is a general fact of commutative algebra.\acapo
We now prove $c)\solose a)$. 
If $r=1$ there is nothing to prove; we begin 
by considering the case $r=2$.\acapo
Every fibre is defined in $\com^3_w$ by the six equations 
$$\cases{w_i^2=f_jf_k&\cr w_jw_k=w_if_i&\cr}\qquad 
\{i,j,k\}=\{1,2,3\},\qquad f_i\in \com$$
If $X$ is not locally simple at a point $x$, then the fibre containing 
$x$ has $f_1=f_2=f_3=0$ and then it is isomorphic to 
$Spec(\com[w_1,w_2,w_3]/(w_1,w_2,w_3)^2)$ 
and then it is not Gorenstein because the 
socle has dimension 3.\acapo
As the property of being Gorenstein is local the same conclusion holds if 
$r\ge 2$ but the stabilizer of every point $x\in X$ is a subspace of 
dimension $\le 2$.\acapo  
In general assume $r\ge 3$ and  $\pi$ not locally simple at $x$; it is 
not restrictive to shrink  $Y$ to a sufficiently small affine 
neighbourhood $U$ of $y=\pi(x)$ such that $D_\sigma\cap U$ are principal 
divisors.\acapo
Let $\sigma,\tau\in G$ such that $\sigma,\tau,\sigma+\tau\in I_x$.
Since the property of being Gorenstein is stable under small deformations, 
we can deform 
the divisors $D_\sigma$ by moving away from $y$ the divisors $D_\alpha$ 
with $\alpha\not=\sigma,\tau,\sigma+\tau$. After this deformation 
$I_x=\{\sigma,\tau,\sigma+\tau\}$ the stabilizer of $x$ has dimension 2 
and the above computation gives a contradiction.\finedim  

\Cor{3.17}{} {\it Let $\pi\colon X\to Y$ be a $\zduer$-cover of algebraic 
surfaces. If $Y$ is smooth and $X$ has at most rational double points as 
singularities then $\pi$ is locally simple.}\acapo
\Proof Immediate from Proposition 3.12.\finedim 

\Ex{3.18} Let $(X_0,0)$ be the cyclic singularity of type $\ds{1\over 
4}(1,1)$, it is well known [Rie] that the semiuniversal deformation of 
$(X_0,0)$ contains two irreducible components: the Artin component and 
the $\raz$-Gorenstein component.\acapo 
The $\raz$-Gorenstein component is the one-parameter smoothing 
$\pi\colon\X=\Y/\mu_2\to \com_t$ where 
$\Y\subset \com^3\times\com_t$ is defined by the 
equation $uv=t+y^2$ and the group $\mu_2$ is generated by the involution 
$(u,v,y)\to (-u,-v,-y)$.\acapo 
The group $G=\zdue{}\times\zdue{}$ generated by the involutions 
$$\sigma(u,v,y,t)=(v,u,y,t),\qquad \varrho(u,v,y,t)=(u,v,-y,t)$$
acts on $\X$ preserving fibres. 
A set of generator of the $\com$-algebra $\O_{\Y}^{\mu_2}=\O_{\X}$ which 
are eigenvectors for the $G$-action is given by 
$$t,\quad x=uv=y^2+t,\quad z=u^2+v^2,\quad w_1=uy+vy,\quad 
w_2=u^2-v^2,\quad w_3=uy-vy$$
It is trivial to see that $\X/G$ is smooth and $\X$ is isomorphic to the 
subvariety of $\com^5\times\com_t$ defined by the equations
$$w_i^2=f_jf_k,\qquad w_iw_j=w_kf_k,\qquad \{i,j,k\}=\{1,2,3\}$$
where $f_1=z-2x$, $f_2=x+t$, $f_3=z+2x$.\acapo
It is important for us to observe that the singularity $X_0$ is the 
bidouble cover of $\com^2$ branched over three lines passing through 0 
and that every $\raz$-Gorenstein smoothing is obtained as a Galois 
deformation by moving a line away from 0.\acapo
\Ex{3.19} Let $m=(C_1,...,C_k,p_1,...,p_n)\in\bar{M}_{k,n}$ and denote by
$\nu\colon S_m\to\pro^1\times\pro^1$ the blow up of $p_1,...,p_n$, 
$E_i=\nu^{-1}(p_i)$ and $D_i=\nu^*C_i-\sum_j E_j$.\acapo
There is a natural morphism of functors of Artin rings 
$\phi\colon (\bar{M}_{k,n},m)\to Def_{S_m,D_i}$ 
and by Kodaira stability theorem $\phi$ is a smooth morphism.\acapo
\Ex{3.20} In the notation of example 3.19 assume $k=3$ and let 
$\pi\colon X_m\to S_m$ be the 
$\zdue{2}=\{0,\alpha_1,\alpha_2,\alpha_3\}$-cover 
with building data 
$D_{\alpha_i}=D_i$, $i=1,2,3$ and $L_{\chi}=\O_{S_m}(D_1)$, 
$\chi\not=0$.\acapo 
In view of  3.7,  3.17 and 3.18 
we see easily that:\acapo 
\item{i)} if $m\in M^0_{3,n}$ then $X_m$ is smooth.\acapo
\item{ii)} if $m\in M_{3,n}$ then $X_m$ has at most cyclic singularities 
of type $\ds{1\over 4}(1,1)$.\acapo
\item{iii)} if
$X$ has at most rational double points as singularities then
$m\in\bar{M}^0_{3,n}$.\acapo 
This partially motivates the definitions of section 2. It is an 
interesting question to decide whether two surfaces $X_m$, $X_l$ are 
deformation equivalent whenever $m,l$ belongs to different connected 
components of $M^0_{3,n}$. We suspect that the answer is no in general 
but a proof of this fact seems inaccessible at the moment.\acapo

\Parag{4}{$\zduer$-actions on rational double points and their smoothings}
In this section we continue the investigation of [Ca3], [Ma3] concerning 
the quotients of rational double points by commuting involutions: in 
particular we are interested to simplicity and smoothability of the 
actions.\acapo
We first introduce some terminology. 
Assume the group 
$G=\zduer$ acts faithfully on a normal surface singularity 
$(X,x)$ and let $\pi\colon (X,x)\to (Y,y)$ be the quotient map.\acapo 
For every $\sigma\in G$ let $R_\sigma\subset X$ be the germ of the 
divisorial part of the fixed locus of $\sigma$. As in the global case, 
every irreducible curve germ through $X$ may belong to at most one 
$R_\sigma$. We denote by $D_\sigma$ the image of $R_\sigma$ with the 
reduced 
structure and $I_x=\{\sigma|R_{\sigma}\not=0\}$.\acapo
\Def{4.1} The $G$-action is called {\it almost simple} if $I_x$ is a set 
of linearly independent vectors. It is called {\it simple} if it is 
almost simple and the quotient map $\pi$ is flat.\acapo
The computations of \S 3 show that if the $G$-action is simple then the 
divisors $R_\sigma$, $D_\sigma$ are principal and $I_x$ is a basis. In 
particular if $Y$ is smooth then the $G$-action is simple if and only if 
it is almost simple.\acapo
After [Ca3] we know that there exist exactly 13 conjugacy classes of 
$\zduer$-actions on rational double points. We shall give this list in 
the next three tables.\medskip\acapo
{\bf Table 4.2}. Equations of RDP's in $\com^3$.\acapo
$$\vbox{\settabs 3\columns
\+ $E_8$&$ z^2+x^3+y^5=0$&\quad\hfil\cr
\+$E_7$&$ z^2+x(y^3+x^2)=0$&\hfil\cr
\+$E_6$&$z^2+x^3+y^4=0$&\hfil\cr
\+$D_n,\, n\ge 3$&$ z^2+x(y^2+x^{n-2})=0$&\hfil\cr
\+$A_n,\, n\ge 0$&$ z^2+x^2+y^{n+1}=0\quad\hbox{or}\quad 
uv+y^{n+1}=0$&\hfil\cr}$$
By a little abuse of terminology, we consider the smooth germ $(\com^2,0)$ 
as the 
rational double point $A_0$ and $D_3=A_3$; clearly $x,z$ are local 
coordinates on 
$\com^2$ at 0. Moreover we have $u=z+ix$, $v=z-ix$.\medskip\acapo

\acapo
{\bf Table 4.3.} ([Ca3] ) The six canonical forms 
of involutions acting 
on the RDP's.\acapo
$$\vbox{\settabs 3\columns
\+$a)$&$(x,y,z)\to (x,-y,z)$&\quad\cr
\+$b)$&$(x,y,z)\to (x,-y,-z)$&\hfil\cr
\+$c)$&$(u,v,y)\to (-u,v,-y)$&\hfil\cr
\+$d)$&$(x,y,z)\to (-x,y,-z)$&\hfil\cr
\+$e)$&$(u,v,y)\to (-u,-v,-y)$&\hfil\cr
\+$f)$&$(x,y,z)\to (x,y,-z)$&\hfil\cr}$$
\medskip\acapo
{\bf Table 4.4.} Conjugacy classes of $\zduer$-actions on rational double 
points.\nobreak 
$$\vbox{\offinterlineskip\hrule\halign
{&\vrule#&\strut~#\hfil~\cr
height4pt&\omit&&\omit&&\omit&&\omit&&\omit&&\omit&&\omit&&\omit&\cr
&&&$r$&&$X$&&basis&&$Y$&&$|I_x|$&&simple?&&smoothable?\hfil&\cr
height4pt&\omit&&\omit&&\omit&&\omit&&\omit&&\omit&&\omit&&\omit&\cr
\noalign{\hrule}
height4pt&\omit&&\omit&&\omit&&\omit&&\omit&&\omit&&\omit&&\omit&\cr
&1 &&$1$ &&$E_6,D_n,A_{2n+1}$&&$a$&&$A_2,A_1,A_{n+1}$&&$1$&&\hfil 
yes&&\hfil yes&\cr 
height4pt&\omit&&\omit&&\omit&&\omit&&\omit&&\omit&&\omit&&\omit&\cr
\noalign{\hrule}
height4pt&\omit&&\omit&&\omit&&\omit&&\omit&&\omit&&\omit&&\omit&\cr
&2 &&$1$&&$E_6,D_n,A_{2n+1}$&&$b$&&$E_7,D_{2n-2},D_{n+3}$&&$0$&&\hfil 
no&&\hfil no&\cr 
height4pt&\omit&&\omit&&\omit&&\omit&&\omit&&\omit&&\omit&&\omit&\cr
\noalign{\hrule}
height4pt&\omit&&\omit&&\omit&&\omit&&\omit&&\omit&&\omit&&\omit&\cr
&3 &&$1$ &&$A_{2n}$  &&$ c$ && $B_n$&& $1$ &&\hfil no&&\hfil yes&\cr 
height4pt&\omit&&\omit&&\omit&&\omit&&\omit&&\omit&&\omit&&\omit&\cr
\noalign{\hrule}
height4pt&\omit&&\omit&&\omit&&\omit&&\omit&&\omit&&\omit&&\omit&\cr
&4 &&$1$ &&$A_n$  &&$ d$ && $A_{2n+1}$&& $0$ &&\hfil no&&\hfil no&\cr 
height4pt&\omit&&\omit&&\omit&&\omit&&\omit&&\omit&&\omit&&\omit&\cr
\noalign{\hrule}
height4pt&\omit&&\omit&&\omit&&\omit&&\omit&&\omit&&\omit&&\omit&\cr
&5 &&$1$ &&$A_{2n-1}$  &&$ e$ && $Y_{n}$&& $0$ &&\hfil no&&\hfil yes&\cr 
height4pt&\omit&&\omit&&\omit&&\omit&&\omit&&\omit&&\omit&&\omit&\cr
\noalign{\hrule}
height4pt&\omit&&\omit&&\omit&&\omit&&\omit&&\omit&&\omit&&\omit&\cr
&6 &&$1$ &&all RDP's &&$ f$ && $A_0$&& $1$ &&\hfil yes&&\hfil yes&\cr 
height4pt&\omit&&\omit&&\omit&&\omit&&\omit&&\omit&&\omit&&\omit&\cr
\noalign{\hrule}
height4pt&\omit&&\omit&&\omit&&\omit&&\omit&&\omit&&\omit&&\omit&\cr
&7 &&$2$ &&$E_6,D_n,A_{2n+1}$  &&$a,f $ && $A_0$&& $2$ &&\hfil yes&&\hfil 
yes&\cr 
height4pt&\omit&&\omit&&\omit&&\omit&&\omit&&\omit&&\omit&&\omit&\cr
\noalign{\hrule}
height4pt&\omit&&\omit&&\omit&&\omit&&\omit&&\omit&&\omit&&\omit&\cr
&8 &&$2$ &&$A_n$  &&$ d,f$ && $A_0$&& $2$ &&\hfil yes&&\hfil yes&\cr 
height4pt&\omit&&\omit&&\omit&&\omit&&\omit&&\omit&&\omit&&\omit&\cr
\noalign{\hrule}
height4pt&\omit&&\omit&&\omit&&\omit&&\omit&&\omit&&\omit&&\omit&\cr
&9 &&$2$ &&$A_{2n+1}$  &&$a,d$ && $A_{2n+1}$&& $1$ &&\hfil no&&\hfil 
no&\cr 
height4pt&\omit&&\omit&&\omit&&\omit&&\omit&&\omit&&\omit&&\omit&\cr
\noalign{\hrule}
height4pt&\omit&&\omit&&\omit&&\omit&&\omit&&\omit&&\omit&&\omit&\cr
&10 &&$2$ &&$A_n$  &&$ e,f$ && $A_{1}$&& $1$ &&\hfil no&&\hfil no&\cr 
height4pt&\omit&&\omit&&\omit&&\omit&&\omit&&\omit&&\omit&&\omit&\cr
\noalign{\hrule}
height4pt&\omit&&\omit&&\omit&&\omit&&\omit&&\omit&&\omit&&\omit&\cr
&11 &&$2$ &&$A_{2n+1}$  &&$ b,d$ && $D_{2n+4}$&& $0$ &&\hfil no&&\hfil 
no&\cr 
height4pt&\omit&&\omit&&\omit&&\omit&&\omit&&\omit&&\omit&&\omit&\cr
\noalign{\hrule}
height4pt&\omit&&\omit&&\omit&&\omit&&\omit&&\omit&&\omit&&\omit&\cr
&12 &&$2$ &&$A_{2n+2}$  &&$c, d$ && $A_{2n+2}$&& $2$ &&\hfil no&&\hfil 
yes&\cr 
height4pt&\omit&&\omit&&\omit&&\omit&&\omit&&\omit&&\omit&&\omit&\cr
\noalign{\hrule}
height4pt&\omit&&\omit&&\omit&&\omit&&\omit&&\omit&&\omit&&\omit&\cr
&13 &&$3$ &&$A_{2n+1}$  &&$a,d,f$ && $A_0$&& $3$ &&\hfil yes&&\hfil 
yes&\cr 
height4pt&\omit&&\omit&&\omit&&\omit&&\omit&&\omit&&\omit&&\omit&\cr}
\hrule}$$\acapo
The above Table 4.4 requires some explanations.\acapo
By definition $B_n$ is the quotient of the rational 
double point $A_{2n}$ by the involution $c)$; the computation of [Ca3] 
shows that $B_n$ is the cyclic  singularity of type 
$\ds{1\over 2n+1}(1,2n-1)$ and by the Hirzebruch-Jung expansion as a 
continued fraction ([B-P-V] III.5) its Dynkin diagram is 
$$\cecc{-2}\riga\cecc{-2}\riga\cecc{-2}
\ldotsa\cecc{-2}\riga\cecc{-3}\quad (n\ge 2\,\hbox{vertices})\eqno{(4.5)}$$
If $Z$ is the fundamental cycle then $Z^2=-3$ and then $B_n$ is a 
rational triple point; by Riemenschneider results ([Rie, Satz 10 and Satz 
12]) 
every deformation 
of $B_n$ admits simultaneous resolution. In particular if $F$ is the 
Milnor fibre of a smoothing of $B_n$ then $H_2(F,\interi)$ is torsion free 
of rank $n$ with the intersection form induced by the Dynkin diagram 
$(4.5)$.\acapo 
By definition $Y_{n}$ is the quotient of the 
rational double point $A_{2n-1}$ by the involution $e)$. A simple 
computation ([Ca3, p. 92], [Ma1], [Ma7]) shows that $Y_{n+1}$ is the 
cyclic quotient 
singularity of type $\ds{1\over 4n}(1,2n-1)$ and its Dynkin diagram is 
$$\cecc{-4}\qquad\hbox{ for $n=1$}$$
$$\cecc{-3}\riga\cecc{-2}\riga\cecc{-2}
\ldotsa\cecc{-2}\riga\cecc{-3}\quad (n\ge 2\,\hbox{vertices})\qquad 
\hbox{ for $n\ge 2$}$$
The fundamental cycle has selfintersection $Z^2=-4$ and then $Y_{n}$ is a 
rational quadruple point. Again by Riemenschneider results ([Rie, Satz 
13]) the base space of the semiuniversal deformations of $Y_{n}$ is the 
transversal union of two smooth germs parametrizing respectively 
deformations admitting simultaneous resolution and $\raz$-Gorenstein 
deformations.\acapo 
If $F$ is the Milnor fibre of a $\raz$-Gorenstein smoothing of $Y_n$ then 
it is proved in ([Ma1, Prop. 13]) that the torsion subgroup 
of $H^2(F,\interi)$ is isomorphic to $\zdue{}$ and it is generated by the 
canonical class.\acapo     
The first 4 columns of 4.4, namely $r$, 
$X$, $Y$ and the basis of $G$; come directly from [Ca3]. By $|I_x|$ we 
denote the cardinality of the set $I_x$: a case by case computation 
shows that all the above actions are almost simple.\acapo
The sixth column tell us if the action is simple, this is obtained by  
direct inspection. The last column tell us if the action is 
smoothable.\acapo
For reader convenience we recall the notion of smoothable action 
introduced in [Ma3]. 
\Def{4.6} Let $G$ be a finite group acting on a normal twodimensional 
singularity $(X,0)$. The action is called {\it smoothable} if there 
exists a smoothing $(\X,0)\to(\com,0)$ of $(X,0)$ such that the 
$G$-action extends to $\X$ preserving fibres  and the quotient 
$(\X/G,0)\to (\com,0)$ is a smoothing of the quotient $Y=X/G$. The Milnor 
fibre of $\X/G\to\com$ is called an {\it admissible Milnor fibre} of the 
$G$ action.\acapo
\Rem{4.7} Not every Milnor fibre of singularity $X/G$ is admissible, 
for example if $(X,0)$ is a rational double point and $G$ is a cyclic 
group 
acting freely on $X-0$ then the admissible Milnor fibres are exactly the 
ones coming from $\raz$-Gorenstein smoothing.\acapo
In the definition 4.6, the proof that $Y=X/G$ is the central fibre of the 
smoothing $\X/G$ is non trivial and rests on the Serre criterion 
[Mat, 23.8] and the 
fact that Cohen-Macaulay complex spaces are preserved under finite group 
actions.\acapo 
The answer yes in the last column of 4.4 
is given by writing down explicit smoothings, 
while the answer no in obtained by using the same methods of [Ma3] 
(see below).\acapo
Note that every simple action realize $X$ as a simple cover of $Y$ and 
then a generic Galois deformation of the cover gives a smoothing of the 
action (recall that every rational twodimensional singularity is 
smoothable,
[Wa3, 1.3]).\acapo 
In the following example we describe some smoothings of non simple 
actions (i.e. lines 3,5 and 12 of Table 3).
\Ex{4.8} {} The singularity $\X\subset \com^3\times\com$ defined by the 
equation $uv+y^{2n+1}+ty=0$ defines a smoothing of the rational double 
point $A_{2n}$; the involution $c)$ extends to $\X$ and its quotient is a 
smoothing of $B_n$.\acapo 
The negative answer to smothability in lines 2,4 is proved in 
[Ma3, 4.6], the negative answer in lines 9,10,11 is proved in a similar 
way; 
here we illustrate the idea by showing that the action 9 is not 
smoothable.\acapo 
The rational double point $X$ of type $A_{2n-1}$, $n>0$, is defined in 
$\com^3$ by the equation $z^2+x^2+y^{2n}=0$ and the action 9 is given by 
the involutions $a,d,e=ad$, $G=\{1,a,d,e\}$.\acapo
Since $G$ is finite, every $G$-stable smoothing of $X$ over $(\com,0)$ is 
defined in $(\com^3\times\com,0)$ by an equation 
$$F_t(x,y,z)=z^2+x^2+y^{2n}+t\phi(x,y,z,t)=0$$
with $\phi$ a $G$-invariant function; see [Ma3, 4.5] for a proof of 
this fact.\acapo
In particular the power series expansion of $\phi$ 
cannot contain linear terms in $x,y,z$. Since for $t\not=0$ sufficiently 
small the fibre $X_t=\{F_t=0\}$ is smooth, we must have necessarily 
$\phi(0,0,0,t)\not=0$ and the intersection of $\X$ 
with the line $\{x=z=0\}=Fix(d)$ is the germs of plane curve of equation 
$y^{2n}+t\phi(0,y,0,t)=0$. In particular, for $0<|t|<<1$ there exists a 
point $x=(0,\xi,0,t)\in Fix(d)\cap X_t$ with $\xi\not=0$; hence $x$   
is an isolated fixed point for the involution 
$d$ and $x\not\in Fix(a)\cup Fix(e)$. 
This implies that the image of 
$x$ in the quotient $X_t/G$ is a singular point.\acapo
In the following three lemmas we deduce, in the style of [Ca3], some 
topological obstructions to degenerations.\acapo

\Lemma{4.9}{} {\it Let $V\mapor{}\com^2$ be the blow-up at a point and let 
$F_n$ 
be the Milnor fibre of a smoothing of a singularity of type 
$A_n,D_n,E_n,B_n$, $n>0$.\acapo
Then cannot exist any open embedding $F\subset V$.}\acapo

\Proof    Assume that $F_n$ is an open subset of $V$, then there is a 
homomorphism of groups $i_*\colon H_2(F_n,\interi)\to H_2(V,\interi)$ 
which preserves 
the intersection products.  We have $H_2(V,\interi)=\interi e$, where $e$ 
is the class of the exceptional curve, $e^2=-1$. By simultaneous 
resolution the intersection form on $H_2(F_n,\interi)$ is negative 
definite 
of rank $n$; this forces $n=1$.
But $H_2(F_1,\interi)=\interi \epsilon$ with $\epsilon^2=-2$ for 
smoothing of $A_1$ and $\epsilon^2=-3$ for smoothing of $B_1$; 
in both cases $i_*$ cannot be an isometry.\finedim

\Lemma{4.10}{} {\it Let $F$ be the Milnor fibre of a smoothing of $Y_n$, 
then 
cannot exist any open embedding $F\subset \com^2$.}\acapo
\Proof There are two different cases:\acapo
1) $F$ is the fibre of a smoothing admitting simultaneous resolution. In 
this case the intersection form on $H_2(F,\interi)$ is negative definite 
and cannot exist any isometry $H_2(F,\interi)\to 
H_2(\com^2,\interi)=0$.\acapo
2) $F$ is the fibre of a $\raz$-Gorenstein smoothing. As $F$ is open in 
$\com^2$ the (trivial) canonical class restricts to the canonical class 
of $F$ which is non trivial; this gives a contradiction.\finedim

\Lemma{4.11}{} {\it Let $F$ be the Milnor fibre of a smoothing of a 
rational 
double point of type $A_n$, $n\ge 2$ and let $V\to \com\times\pro^1$ be 
the blow up at two points. Then cannot exist any open embedding 
$F\subset V$.}\acapo 
\Proof As above an embedding would give an isometry     
$i_*\colon H_2(F,\interi)\to H_2(V,\interi)=\interi^3$. Since the 
intersection form on $H_2(V,\interi)$ is given by the diagonal matrix 
$diag(0,-1,-1)$ a simple calculation shows that this leads to a 
contradiction.\finedim

\Parag{5}{A global construction} 
Let $k>0$ be a fixed natural number  
and consider integers $a_1,...,a_k,b_1,...,b_k\ge 3$,
$n_1<n_2<...<n_k$, $a_i\not=b_i$ such that $0\le n_i\le 2a_ib_i$
and denote $\bar{M}_{3,n_i}=\bar{M}_{3,n_i}^{a_i,b_i}$.\acapo
Let $\bar{M}\subset \bar{M}_{3,n_1}^{a_1,b_1}\times...\times
\bar{M}_{3,n_k}^{a_k,b_k}$ be 
the open subscheme of $k$-uples $((C^1_i,p^1_j),...,(C^k_i,p^k_j))$ such 
that $p^i_j\not=p^h_l$ for $(i,j)\not=(h,l)$ and the curves $C^i_j$ 
are without common components.\acapo
Define also 
$$\bar{M}^0=\bar{M}\cap\prod_i\bar{M}_{3,n_i}^0.$$
Let $\nu\colon \SS\to \bar{M}\times Q$ be the blow up of the 
sections $p^i_j\colon \bar{M}\to \bar{M}\times Q$, $i=1,...,k$ 
$j=1,...,n_i$ 
and denote by $p\colon \SS\to \bar{M}$, 
$q=(q_1,q_2)\colon\SS\to \pro^1\times\pro^1=Q$ 
the natural projections.\acapo
Denote by $E^i_j\subset \SS$  the exceptional Cartier divisor over the 
sections $p^i_j$. We have seen that the coherent sheaf 
$p_*(q^*\O_Q(r,s)-\sum a^i_jE^i_j)$ is locally free in the following 
cases:\acapo
\item{i)} $r,s\ge 0$ and $a^i_j=0$ for every $i,j$.\acapo
\item{ii)} $a^i_j\ge 2$ for every $i,j$ and $r,s>\sum(a^i_j+1)$.\acapo
Moreover in the above cases i), ii) the restriction of $q^*\O_Q(r,s)-\sum 
a^i_jE^i_j$ to every fibre of $p$ gives a base point free linear 
system.\acapo
For every $i=1,...,k$ define
$G^i_1=\zdue{2}=\{0,\alpha^i_1,\alpha^i_2,\alpha^i_3\}$, 
$G^i_2=\zdue{n_i}$ 
with basis $\epsilon^i_1,...,\epsilon^i_{n_i}$. Define moreover 
$G\pr=G^1_1\oplus G^1_2\oplus....\oplus G^k_1\oplus G^k_2\oplus 
\zdue{4}$ with $\tau_1,\tau_2,\eta_1,\eta_2$ a basis of $\zdue{4}$ 
and $G=G\pr\oplus\zdue{}$.\acapo
Given a point $m=((C^1_i,p^1_j),...,(C^k_i,p^k_j))\in \bar{M}$ 
we want to construct a 
family of normal flat $G$-covers over the fibre  
$S_m=p^{-1}(m)$, 
the building 
data $L\colon G\dual\to Pic(S_m)$, $D\colon G\to Div(S_m)$
of this cover must satisfy the following conditions:\acapo
\item{i)} $L_\chi-D_\sigma$ is combinatorially ample for every 
$\chi\not=0$, $\sigma\in G$.\acapo
\item{ii)} $D_{\alpha^i_j}=q^*C^i_j-\sum_h E^i_h$ $j=1,2,3$, 
$i=1,..,k$.\acapo 
\item{iii)} $D_{\epsilon^i_j}=E^i_j$, $D_{\tau_i}$ a fibre of $q_1$, 
$D_{\eta_i}$ a fibre of $q_2$. If $\sigma\in G\pr$ is not one 
of the previous cases then $D_\sigma=0$.\acapo
\item{iv)} If $\sigma\in G-G\pr$ then $D_\sigma$ is 
combinatorially ample.\acapo
We have seen that there exists sufficiently ample building data 
satisfying the above conditions, moreover by the computation of section 2 
it follows easily that the triple $(S,L,D)$ satisfies the hypothesis of 
corollary 3.14.\acapo
Since $Pic(S_m)$ is torsion free the map $L\colon G\dual\to Pic(S_m)$ is 
uniquely determined by $D\colon G\to Div(S_m)$. For a fixed $m$ and $L$ 
the 
set of building data $(L,D)$ as above is then parametrized by a product 
$U_m$ of projective spaces; by base change the union of all $U_m$, with 
fixed $L$ is an algebraic  bundle $s\colon U\to \bar{M}$.\acapo
\Lemma{5.1}{} {\it Let $U^0\subset U$ (resp.: $U^{00}$) 
be the subset of triples $(m,D,L)$ 
such that the divisors $D_\sigma$ are reduced without common components 
and the associated $G$-cover $X\to S_m$ has at most rational double 
points as singularities (resp.: $X$ is smooth).\acapo 
Then $U^0, U^{00}$ are Zariski open in $U$.\acapo
Moreover $ s(U^{00})\subset s(U^0)\subset \bar{M}^0$.}\acapo
\Proof For every $u=(s(u),L_u,D_u)\in U$ 
denote $S_u=\SS\times_{\bar{M}}\{u\}$. 
By construction there exist (tautological) maps $D_U\colon G\to 
Div(\SS\times_{\bar{M}}U)$, 
$L_U\colon G\dual\to Pic(\SS\times_{\bar{M}}U)$ such that their 
restriction to $S_u$ are exactly $D_u,L_u$; define $\tilde{U}\subset U$ 
as the Zariski open subset consisting of points $u$ such that the 
divisors $D_{u,\sigma}$ are reduced without common components and then 
they are the branching divisors of a normal flat $G$-cover $X_u\to 
S_u$.\acapo
There exists a Zariski open covering $\tilde{U}=\cup V_i$ such that the 
maps $D_U,L_U$ satisfies the cover conditions on 
$\SS\times_{\bar{M}}V_i$ and therefore there exist flat $G$-covers 
$X_{V_i}\to \SS\times_{\bar{M}}V_i$, which by construction contains 
all the Galois deformations of $X_u\to S_u$ for every $u\in 
\tilde{U}$.\acapo 
By the stability property of the classes of smooth and  rational 
double points under 
deformations it follows that 
$V_i\cap U^0$ and $V_i\cap U^{00}$   
are  Zariski open subsets of $V_i$ and this 
proves that $U^0, U^{00}$ are Zariski open in $U$.\acapo
For $u\in U^0$ let $X_u\to S_u=\SS\times_{\bar{M}}\{u\}$ be the 
associated $G$-cover, by Corollary 3.17 it is locally simple, hence for 
every $i=1,...,k$,  
$D_{\alpha^i_1}\cap D_{\alpha^i_2}\cap D_{\alpha^i_3}=\vuoto$ and this is 
equivalent to the fact that $s(u)\in \bar{M}^0$.\finedim 
If $\M$ denotes the moduli space of surfaces of general type the above 
proof shows that the natural map $\phi\colon U^0\to\M$, $\phi(u)=[X_u]$ 
is a regular morphism of varieties.
Note that the map $\phi$ is invariant under the natural action of 
$Aut_0(Q)$ over $U^0$; it is also useful to remark that $U^0$ and 
$U^{00}$ are stable under the action of 
$Aut_0(Q)\times \prod_{i=1}^k (\Sigma_3\times \Sigma_{n_i})$
but in general $\phi$ is not invariant under this extended action.\acapo 
\Prop{5.2}{} {\it The morphism $\phi$ is open.}\acapo
\Proof Immediate from the fact that, as the building data are choosed
satisfying the hypothesis of 3.14, the Galois deformations are 
complete.\finedim 
Next we want to prove that  $\phi(V)$ is closed in $\M$
for every irreducible component $V$ of $U^0$; since 
$U^{00}$ is dense in $U^0$ it is sufficient to prove that there exists an 
open covering $U^{00}=\cup V_i$ such that $\bar{\phi(V_i\cap V)}\subset 
\phi(V)$ for every $i$.\acapo
We choose the covering $\{V_i\}$ such that there exist global 
$G$-covers $\X_{V_i}\to \SS\times_{\bar{M}}V_i$, by the valuative 
criterion the closure of $\phi(V)$ follows from the following 
\Theor{5.3}{} {\it Let $U^{00}=\cup V_i$ be as above and 
let $\Delta$ be a smooth affine curve, $0\in \Delta$ 
and let $f\colon \X\to \Delta$ be a proper flat family of irreducible 
surfaces 
such that:\acapo
\item{i)} $X_t=f^{-1}(t)$ is smooth for $t\not=0$ and $X_0$ has at most 
rational double points as singularities.\acapo
\item{ii)} $X_t$ has ample canonical bundle for every $t\in \Delta$.\acapo
\item{iii)} Let $\X^*=\X-X_0$, $\Delta^*=\Delta-\{0\}$, then 
there exists a regular morphism $\eta\colon \Delta^*\to V_i\subset U^{00}$
such that $\X^*\to \Delta^*$ is isomorphic to the pull-back of 
$\X_{V_i}\to V_i$ under $\eta$.\acapo
Then the map $\eta\colon \Delta^*\to U^{00}$ can be extended, up to the 
action of $Aut_0(Q)$ on $U^{00}$,  
to a regular morphism $\eta\colon \Delta\to U^0$.}\acapo 
\Proof  As the fibres of $f$ have at most rational double points and 
ample canonical bundle the hypothesis of the semicontinuity theorem [FP, 
4.4] are satisfied and the $G$-action over $\X^*$ can be extended to a 
regular $G$-action on $\X$; moreover the restriction of this action to 
$X_0$ is faithful.\acapo
Let $\pi\colon\X\to\Y$ be the projection to the quotient and let
$g\colon \Y=\X/G\to\Delta$ be the factorization of $f$. Denote by 
$Y_t=g^{-1}(t)$; note that if $Y_t$ is smooth or $X_t/G$ is smooth and 
$Y_t$ is normal then the natural morphism 
$X_t/G\to Y_t$ is an isomorphism of varieties.\acapo
{\it Claim 1. If $x\in X_t$ is a smooth point then $g\colon \Y\to \Delta$ 
is smooth at $\pi(x)$, in particular $\Y$ is normal with isolated 
singularities.}\acapo 
We can check the smoothness of $g$ from the analytic point of view. Let 
$G_x\subset G$ be the stabilizer of $x$, according to [Ma3, 4.5] there 
exists a $G_x$-isomorphism of analytic germs $(\X,x)=(X_t,x)\times 
(\Delta,t)$ where $G_x$ acts trivially on $\Delta$. 
Then $(\Y,\pi(x))=(X_t/G_x,x)\times (\Delta,t)$, using the fact that 
$Y_t$ is smooth for every $t\not=0$ it follows that $(X_t/G_x,x)$ must be 
smooth. Since $\X$ is normal also $\Y$ is normal.\acapo 
{\it Claim 2. $g\colon \Y\to \Delta$ is smooth and $\pi\colon \X\to \Y$ 
is flat.}\acapo 
By Claim 1 if $g$ is not smooth at a point $y=\pi(x)$ then $x$ must be a 
singular point of $X_0$. As above let $G_x$ be the stabilizer of $x$; 
then $G_x$ acts on the smoothing $(\X,x)\to (\Delta,0)$ preserving fibres 
and its quotient is the smoothing $(\Y,y)\to (\Delta,0)$. Therefore the 
$G_x$-action on the rational double point $(X_0,x)$ is smoothable.\acapo 
For every $\sigma\in G$ let $D_\sigma\subset \Y$ be the reduced branching 
divisor of the action of $\sigma$ on $\X$. If $F\subset Y_t$ is the 
Milnor fibre of $(\Y,y)\to (\Delta,0)$ then for $t$ sufficiently near to 
$0$ we have $F\subset Y_t-(\ds\cup_{\sigma\not\in I_x}D_\sigma)$ and 
then, recalling that among the divisors $D_\sigma\cap Y_t$ there are 
fibres of the two projections 
$q_1,q_2\colon Y_t\to \pro^1\times\pro^1\to \pro^1$, we 
see that $F$ is biholomorphic to an open subset of $V$, where 
\item{a)} $V=\com^2$ if $|I_x|=0$, or 
\item{b)} $V=$ the blow up of $\com^2$ at a point if $|I_x|=1$, or 
\item{c)} $V=$ the blow up of $\com\times \pro^1$ at two points if 
$|I_x|=2$.\acapo 
By using table 4.4 and the topological obstructions 4.9, 4.10 and 4.11 we 
deduce that $Y_0$ must be smooth at $y$.\acapo 
By 3.5 the divisors $D_\sigma\cap Y_t$ are reduced without common 
components and the same argument used in Claim 1 shows that 
for every $t$  every  branching divisor 
of $X_t\to Y_t$ has the form $D_\sigma\cap Y_t$. 
Since $D_{\tau_1}\per D_{\tau_2}\per Y_t=0$, 
we have $D_{\tau_1}\cap 
D_{\tau_2}=\vuoto$ and, possibly shrinking
$\Delta$, $D_{\tau_1}$ is linearly equivalent to 
$D_{\tau_2}$ and then they define a morphism $\Y\to \pro^1$. 
The same   
argument applies to the divisors $D_{\eta_1}$, $D_{\eta_2}$ and then we 
have a morphism $\nu\colon\Y\to Q\times \Delta$ such 
that for every $t\not=0$ the restriction $\nu\colon Y_t\to Q$ is     
the blow up at $n$ distinct points.\acapo
{\it Claim 3. $D_{\epsilon^i_j}\cap Y_0$ is a smooth rational curve for 
every $i=1,...k$, $j=1,...,n_i$.\acapo}
Let $(i,j)$ be  fixed pair as above 
and let $E=D_{\epsilon^i_j}\cap Y_0$.\acapo 
As $D_{\epsilon^i_j}$ is reduced irreducible and the local rings of $\Delta$ 
are discrete valuation rings the map $g\colon D_{\epsilon^i_j}\to \Delta$ 
is flat and then $E\subset Y_0$ is a reduced connected 
divisor of arithmetic genus $p_a(E)=0$.
It is then easy to see that if  
$E=\cup _i E_i\subset Y_0$ is the decomposition in irreducible components, 
then 
every irreducible component $E_i$ is smooth rational, 
$E_i\per E_j\le 1$ for every $i\not=j$ and the dual intersection graph is 
a tree; a trivial calculation shows that $E^2=\sum_i(E_i^2+2)-2$.\acapo 
The one-dimensional variety $\nu(D_{\epsilon^i_j})$ is irreducible and 
then 
$\nu$ contracts $E$; by Mumford theorem [Mu1] every component $E_i$ has 
negative selfintersection. 
If $E^2_i\le -2$ for every $i$ then we get a contradiction by using the 
above formula $-1=E^2=\sum_i(E_i^2+2)-2$.
Let $E_0\subset E$ be a component with $E_0^2=-1$, then by Kodaira 
stability theorem [Ko], [Ho], there exist an analytic neighbourhood 
$0\in \Delta\pr\subset\Delta$ and a
smooth closed submanifold 
$\E\subset g^{-1}(\Delta\pr)$ such that $\E\cap Y_0=E_0$ 
and $g\colon\E\to \Delta\pr$ is a $\pro^1$ bundle. As $\nu$ contracts 
$E_0$ then it contracts all the fibres of $g\colon\E\to \Delta\pr$ and 
therefore $\E\subset D_{\epsilon^h_l}$ for some $h,l$. 
Using the fact that the divisors 
$D_{\epsilon^h_l}\cap Y_0$ $j=1,...,n$ have no common components we get 
$\E\subset D_{\epsilon^i_j}$; both divisors are irreducible and then 
$\E=D_{\epsilon^i_j}$.\acapo 
From Claim 3 it follows that $\nu\colon Y_0\to Q$ is the blow up at $n$ 
distinct points $p_1,...,p_n$. Moreover $D_{\alpha^i_j}\cap Y_0$ is a 
reduced effective divisor linearly equivalent to $\O_{Y_0}(a_i,b_i)-\sum 
E^i_j$ and $D_{\alpha^i_1}\cap D_{\alpha^i_2}\cap 
D_{\alpha^i_3}=\vuoto$; this clearly implies that  
$$(..,(\nu(D_{\alpha^i_1}\cap Y_0), \nu(D_{\alpha^i_2}\cap 
Y_0),\nu(D_{\alpha^i_3}\cap Y_0), p^i_1,...,p^i_{n_i}),...)\in \bar{M}^0$$
From this we get
the required morphism $\eta\colon \Delta\to U^0$.\finedim

\Cor{5.4}{} {\it 
Let $V$ be an irreducible (resp.: connected) component of $U^0$. Then 
$\phi(V)$ is an irreducible (resp.: connected) component of the moduli 
space $\M$.}\acapo
\Proof 
It is an immediate consequence of the fact that $\phi$ is an open map and 
that $\phi(V)$ is closed in $\M$.\finedim
There exists a natural diagonal action of $Aut_0(Q)^k$ on the spaces 
$\prod M_{3,n_i}$ and $\prod M_{3,n_i}^0$. 
\Lemma{5.5}{} {\it 
Let $U_T\subset \tilde{U}$ be the open subset 
of points $u$ such that $X_u$ is a surface of class T.\acapo 
Then the open subsets $s(U_T)\cap \prod M_{3,n_i}\subset \prod M_{3,n_i}$, 
$s(U^{00})\cap \prod M_{3,n_i}\subset \prod M_{3,n_i}^0$
intersect every $Aut_0(Q)^k$-orbit. In particular 
$s(U_T)\cap \prod M_{3,n_i}$ is connected.}\acapo
\Proof Let $m=(..,(C^i_j,p^i_l),..)\in \prod M_{3,n_i}$, for a generic 
choice of $\phi=(\phi_1,...,\phi_k)\in Aut_0(Q)^k$, the curves 
$$\phi_1(C^1_1+C^1_2+C^1_3),...,\phi_k(C^k_1+C^k_2+C^k_3)$$
intersects transversally. This implies that there exists  $u\in U_T$ such 
that $s(u)=\phi(m)$ and every point of $X_u$ is either smooth or cyclic 
singular of type $\ds{1\over 4}(1,1)$. The same argument shows that if 
$m\in \prod M_{3,n_i}^0$ then for generic $u\in U_{\phi(m)}$ 
the surface $X_u$ is smooth.\acapo
Since $\prod M_{3,n_i}$ is connected and $Aut_0(Q)^k$ is irreducible it is 
an easy exercise to deduce that $s(U_T)\cap \prod M_{3,n_i}$ is 
connected.\finedim
As a corollary of 5.5 we get the following 
\Prop{5.6}{} {\it If $u,v\in U^{00}$ then $X_u$, $X_v$ are deformation 
$T$-equivalent and hence diffeomorphic.}\acapo
\Proof It is an immediate consequence of 5.5 and the inclusion 
$s(U^{00})\subset s(U_T)\cap \prod M_{3,n_i}$.\finedim
It remains to estimate from below the number of connected components of 
$\phi(U^0)$.\acapo 
\Prop{5.7}{} {\it If the building data $(L,D)$ are sufficiently ample 
then there exists an open dense subset $U\pr\subset U^{00}$ such that if 
$u,v\in U\pr$ and $\phi(u)=\phi(v)$ then $s(u)=hs(v)$ for some 
$h\in Aut_0(Q)\times \prod (\Sigma_3\times \Sigma_{n_i})$. In particular 
the number of connected components of $\phi(U^0)$ is greater or equal than
the number of connected components of $\prod M_{3,n_i}^0$ which are 
stable by the action of $Aut_0(Q)\times \prod (\Sigma_3\times 
\Sigma_{n_i})$.}\acapo 
\Proof According to [FP,4.4], if $(L,D)$ are sufficiently ample there 
exists an open dense subset $U\pr\subset U^{00}$ such that for every 
$u\in U\pr$, $Aut(X_u)=G$, i.e. every biregular automorphism of $X_u$ is 
an automorphism of the cover $X_u\to S_u$. (The same result can be 
obtained as a straightforward generalization of the main theorem of 
[Ma4]).\acapo 
The point $\phi(u)\in\M$ determines $X_u$ up to biregular isomorphism; if 
$X\simeq X_u$, $u\in U\pr$ then there exists a commutative diagram with 
horizontal isomorphisms
$$\matrix{X&\mapor{}&X_u\cr \mapver{\pi}&&\mapver{}\cr 
X/Aut(X)&\mapor{f}&S_u\cr}$$
The exceptional curves $f^*(E^i_j)$ are exactly the branching divisors 
$D_\sigma$ of the $G$-cover $\pi$ such that $D^2_\sigma=-1$. By 
considering the branching divisors $D_\alpha$ such that $D_\alpha\per 
f^*(E^i_j)=1$ for some $i,j$, and using the fact that $n_1<n_2<...<n_k$ 
it is easy to recover from the isomorphism class of $X_u$ the orbit of 
$s(u)$ under the action of $Aut_0(Q)\times \prod (\Sigma_3\times 
\Sigma_{n_i})$.\acapo 
This proves the first part of the proposition. The second part follows 
immediately from the fact that $s(U\pr)$ intersects every connected 
component of $\prod M_{3,n_i}^0$ and the fact that if $U_1,U_2\subset U$ 
are connected components then either $\phi(U_1)\cap \phi(U_2)=\vuoto$ or 
$\phi(U_1)=\phi(U_2)$.\finedim   
\Cor{5.8}{} {\it In the above situation for suitable values of 
$a_i,b_i,n_i$ the image $\phi(U^0)\subset \M$ is the union of at least 
$2^k$ irreducible components.}\acapo
\Proof Immediate consequence of 5.7 and 2.5.\finedim
\medskip\acapo
{\ninepoint
{\it Aknoledgements.} It is a pleasure to thank F. Catanese for 
introducing me to this subject (cf. [Ma4]) 
and E. Arbarello for some useful suggestions concerning the use of 
Brill-Noether theory in the study of the moduli space of surfaces.}
\medskip\acapo
\centerline{\bf References}
\medskip\acapo
{\ninepoint
\bib{Ale} V.~Alexeev: {\it Moduli spaces $M_{g,n}(W)$ for surfaces.} (M. 
Andreatta and T. Peternell eds.) Proceedings of the international 
conference, Trento 1994, Walter de Gruyter (1996).\acapo
\bib{Ar} M.~Artin: {\it Deformations of singularities.} Tata Institute 
of fundamental research Bombay (1976).\acapo
\bib{ACGH} E.~Arbarello, M.~Cornalba, P.~Griffiths, J.~Harris:
{\it Geometry of algebraic curves, I.} Springer (1984).\acapo
\bib{BPV} W.~Barth, C.~Peters, A.~van de Ven: {\it Compact complex 
surfaces.} 
Sprin\-ger-Verlag Ergebnisse {\bf 4} (1984).\acapo
\bib{Bon} F.~Bonahon: {\it Diff\'eotopies des espaces lenticulaires.} 
Topology {\bf 22} (1983) 305-314.\acapo
\bib{Car} H.~Cartan: {\it Quotient d'un espace analytique par un groupe 
d'automorphismes.} Algebraic geometry and topology: A symposium in honour 
of S.Lefschetz. Princeton Math. Series {\bf 12} (1957) 90-102.\acapo
\bib{Ca1} F.~Catanese: {\it On the moduli spaces of surfaces of general 
type.} 
J. Diff. Geometry {\bf 19} (1984) 483-515.\acapo
\bib{Ca2} F.~Catanese: {\it Moduli of surfaces of general type.} In: 
{\it Algebraic geometry: open problems. Proc. Ravello 1982} 
Springer L.N.M. {\bf 997} (1983) 90-112.\acapo 
\bib{Ca3} F.~Catanese: {\it Automorphisms of rational double points and 
moduli 
spaces of surfaces of general type.} 
Compositio Math. {\bf 61} (1987) 81-102.\acapo
\bib{Ca4} F.~Catanese: {\it Connected components of moduli spaces.} 
J. Diff. Geometry {\bf 24} (1986) 395-399.\acapo
\bib{Ca5} F.~Catanese: {\it (Some) Old and new results on algebraic 
surfaces.}
Proc. I European congress of Math. Paris 1992, 
Birkhauser (1994).\acapo
\bib{Ca6} F.~Catanese: {\it Moduli of algebraic surfaces.} Springer L.N.M 
{\bf 1337} (1988) 1-83.\acapo
\bib{Ca7} F.~Catanese: {\it Everywhere non reduced moduli space.} Invent. 
Math.
{\bf 98} (1989) 293-310.\acapo
\bib{Ch} M.C.~Chang: {\it The number of components of Hilbert schemes.} 
Int. J. Math. {\bf 7} (1996) 301-306.\acapo
\bib{Do} S.~Donaldson: {\it The Seiberg-Witten equations and 4-manifold 
topology} 
Bull. A.M.S. {\bf 33} (1996) 45-70.\acapo
\bib{EV} H.~Esnault, E.~Viehweg: {\it Two dimensional quotient singularities 
deform to quotient singularities.} 
Math. Ann. {\bf 271} (1985) 439-449.\acapo
\bib{FaMa} B.~Fantechi, M.~Manetti: {\it Obstruction calculus for functors 
of Artin rings, I.} J. Algebra, to appear.\acapo 
\bib{FP} B.~Fantechi, R.~Pardini: {\it Automorphism and moduli spaces of 
varieties with ample canonical class via deformations of abelian covers.}
Comm. in Algebra {\bf 25} (1997) 1413-1441.\acapo
\bib{Fle} H.~Flenner: {\it \"Uber Deformationen holomorpher Abbildungen.} 
Habilitationsschrift, Osnabruck 1978.\acapo
\bib{Fre} M.~Freedman: {\it The topology of four-dimensional manifolds.} 
J. Diff. Geometry {\bf 17} (1982) 357-454.\acapo
\bib{Fri} R.~Friedman: {\it Donaldson and Seiberg-Witten invariants of 
algebraic surfaces.} alg-geom/9605006.\acapo
\bib{FrMo1} R.~Friedman, J.W.~Morgan: {\it Algebraic surfaces and 
four-manifolds: some conjectures and speculations.} 
Bull. Amer. Math. Soc. (N.S.) {\bf 18} (1988) 1-19.\acapo
\bib{FrMo2}R.~Friedman, J.W.~Morgan: {\it Smooth four-manifolds and 
complex 
surfaces.} Ergebnisse der Mathematik {\bf 27} Springer-Verlag
(1994).\acapo
\bib{FS} R.~Fintushel, R.~Stern: {\it Rational Blowdowns of smooth 
4-manifolds.}  Alg-geom/9505018.\acapo 
\bib{Gi} D.~Gieseker: {\it Global moduli for surfaces of general type.}
Invent. Math. {\bf 43} (1977) 233-282.\acapo
\bib{Ho} E.~Horikawa: {\it Deformations of holomorphic maps I.} J. Math.
Soc. Japan {\bf 25} (1973) 372-396; II J. Math. Soc. Japan {\bf 26}
(1974) 647-667; III Math. Ann. {\bf 222} (1976) 275-282.\acapo
\bib{Ko} K.~Kodaira: {\it On stability of compact submanifolds of complex 
manifolds.} Am. J. Math. {\bf 85} (1963) 79-94.\acapo
\bib{KS} J.~Koll\'ar, N.I.~Shepherd-Barron: {\it Threefolds and 
deformations of 
surface singularities.} Invent. Math. {\bf 91} (1988) 299-338.\acapo
\bib{Lo} E.J.N.~Looijenga: 
{\it Isolated singular points on complete intersection.} 
London Math. Soc. L.N.S. {\bf 77} (1983).\acapo
\bib{LW} E.J.N.~Looijenga, J.~Wahl: {\it Quadratic functions and smoothing 
surface 
singularities.} Topo\-logy {\bf 25} (1986) 261-291.\acapo 
\bib{Ma1} M.~Manetti: {\it Normal degenerations of the complex projective 
plane.}
J. Reine Angew. Math. {\bf 419} (1991) 89-118.\acapo
\bib{Ma2} M.~Manetti: {\it On Some Components of Moduli Space of Surfaces 
of General 
Type.} Comp. Math. {\bf 92} (1994) 285-297.\acapo
\bib{Ma3} M.~Manetti: {\it Iterated double covers and connected components 
of 
moduli spaces} Topology {\bf 36} (1996) 745-764.\acapo
\bib{Ma4} M.~Manetti: {\it Degenerations of algebraic surfaces and 
applications 
to moduli problems.} Thesis, Scuola Normale Superiore, Pisa (1995).\acapo
\bib{Ma5} M.~Manetti: {\it Automorphisms of generic cyclic covers.}
Revista Matem\'atica de la Universidad Complutense de Madrid
{\bf 10} (1997) 149-156.\acapo
\bib{Ma6} M.~Manetti: {\it Degenerate double covers of the projective 
plane.}
Preprint Scuola Normale Superiore Pisa (1996), to appear in Proc. 
AGE conference, Warwick 1996.\acapo
\bib{Ma7} M.~Manetti: {\it $\raz$-Gorenstein smoothings  of quotient 
singularities.} 
Preprint Scuola Normale Superiore, Pisa (1990) (unpublished).\acapo
\bib{Mat} H.~Matsumura: {\it Commutative Ring Theory.} Cambridge 
University Press 
(1986)\acapo
\bib{Mi} J.~Milnor: {\it Singular points of complex hypersurfaces.} Ann. 
Math. 
Studies 61, Princeton Univ. Press (1968).\acapo
\bib{Mu1} D.~Mumford: {\it The topology of normal singularities of an 
algebraic 
surface and a criterion for simplicity.} 
Publ. Math. IHES {\bf 9} (1961) 5-22.\acapo
\bib{Mu2} D.~Mumford: {\it Abelian varieties.} Oxford Univ. Press 
(1970).\acapo 
\bib{Par} R.~Pardini: {\it Abelian covers of algebraic varieties.} 
J. Reine Angew. Math. {\bf 417} (1991), 191-213.\acapo
\bib{Rie} O.~Riemenschneider: {\it Deformationen von 
quotientensingularit\"aten 
(nach zyklischen grup\-pen).} Math. Ann. {\bf 209} (1974) 211-248.\acapo
\bib{Sch} M.~Schlessinger: {\it Functors of Artin rings.} Trans. Amer. 
Math. Soc. 
{\bf 130} (1968) 208-222.\acapo
\bib{Ty} A.N.~ Tyurin: {\it Six lectures on four manifolds.} C.I.M.E. 
lectures  ``Transcendental methods in algebraic geometry'' 1994, 
Springer LNM {\bf 1646} (1996).\acapo
\bib{Tya} G.N.~Tyurina: {\it Resolutions of singularities of plane 
deformations of 
double rational points.} Funk. Anal. i. Priloz. {\bf 4} (1970) 77-83.\acapo
\bib{Vie} E.~Viehweg: {\it Quasi-projective moduli for polarized 
manifolds.} Springer-Verlag Ergebnisse Math. Grenz. {\bf 30} (1995).\acapo
\bib{Wa1} J.~Wahl: {\it Smoothings of normal surface singularities.} 
Topology 
{\bf 20} (1981) 219-246.\acapo
\bib{Wa2} J.~Wahl: {\it Elliptic deformations of minimally elliptic 
singularities.} Math. Ann. {\bf 253} (1980) 241-262.\acapo
\bib{Wa3} J.M.~Wahl: {\it Simultaneous resolution of rational 
singularities.}
Comp. Math. {\bf 38} (1979) 43-54.\acapo
}

\medskip\acapo
\settabs\+&xxxxxxxxxxxxxxxxxxxxxxxxxxxxxxxxxxxxxxxxxxxxx&\hfill
&\cr
\+& Scuola Normale Superiore &&\cr
\+& Piazza dei Cavalieri, 7 &&\cr
\+& I 56126 Pisa (Italy) &&\cr
\+& {\tt manetti@cibs.sns.it} &&\cr

\bye